\documentclass[11pt]{article}
\usepackage{amsmath, amsthm, amssymb, eucal, setspace}
\newcommand{\comment}[1]   {{}}
\makeatletter
\renewcommand\section{\@startsection{section}{1}{\z@}%
 						{-3.5ex \@plus -1ex \@minus -.2ex}
						{2ex \@plus.2ex}
						{\large\bfseries}}
\renewcommand\subsection{\setcounter{subsection}{\value{equation}}
						\stepcounter{equation}
					\@startsection{subsection}{2}{\z@}
                          {1.75ex \@plus.5ex \@minus.2ex}%
                           {-.4em}
                           {\textit}}
\def\@seccntformat#1{\@ifundefined{#1@cntformat}%
	{\csname the#1\endcsname\quad} 
	{\csname #1@cntformat\endcsname}} 
\def\section@cntformat{\thesection.~} 
\renewcommand*\l@section{\mdseries\small\@dottedtocline{1}{1.5em}{2em}}
\makeatother
\numberwithin{equation}{section}
\swapnumbers

\newtheorem{theorem}[equation]{Theorem}
\newtheorem{corollary}[equation]{Corollary}
\newtheorem{lemma}[equation]{Lemma}
\newtheorem{proposition}[equation]{Proposition}

\theoremstyle{definition}
\newtheorem{definition}[equation]{Definition}
\theoremstyle{remark}
\newtheorem{remark}[equation]{Remark}
\newtheorem{example}[equation]{Example}

\newcommand{\pz}{(\!(z)\!)}

\newcommand{\GL}{\mathrm{GL}}
\newcommand{\SL}{\mathrm{SL}}
\newcommand{\PSL}{\mathbb{P}\mathrm{SL}}
\newcommand{\Sp}{\mathrm{Sp}}
\newcommand\on{\operatorname}
\newcommand{\Sym}{\mathrm{Sym}}

\newcommand{\cP}{\mathcal{P}}

\newcommand{\GC}{G}
\newcommand{\GR}{G_k}

\newcommand{\cB}{\mathcal{B}}

\newcommand{\cE}{\mathcal{E}}

\newcommand{\cI}{\mathcal{I}}
\newcommand{\cK}{\mathcal{K}}
\newcommand{\cL}{\mathcal{L}}
\newcommand{\cO}{\mathcal{O}}
\newcommand{\cQ}{\mathcal{Q}}
\newcommand{\cR}{\mathcal{R}}
\newcommand{\frb}{\mathfrak{b}}
\newcommand{\frg}{\mathfrak{g}}
\newcommand{\frgR}{\mathfrak{g}_k}
\newcommand{\frn}{\mathfrak{n}}
\newcommand{\frp}{\mathfrak{p}}

\newcommand{\frt}{\mathfrak{t}}
\newcommand{\frz}{\mathfrak{z}}
\newcommand{\frM}{\mathfrak{M}}
\newcommand{\frP}{\mathfrak{P}}

\newcommand{\bC}{\mathbb{C}}
\newcommand{\bD}{\mathbb{D}}
\newcommand{\bE}{\mathbb{E}}

\newcommand{\bP}{\mathbb{P}}
\newcommand{\bQ}{\mathbb{Q}}
\newcommand{\bR}{\mathbb{R}}

\newcommand{\bZ}{\mathbb{Z}}
\newcommand{\gr}{\mathrm{gr}}
\newcommand{\Tr}{\mathrm{Tr}}

\newcommand{\SU}{\mathrm{SU}}
\newcommand{\Ind}{\mathrm{Ind}}	

\newcommand{\mi}{\mathrm{i}}
\newcommand{\waff}{W_\mathrm{aff}}
\newcommand{\hroot}{\vartheta}
\newcommand{\hess}{H}

\begin{document}
\title{\textbf{The index formula for the moduli\\ of $\GC$-bundles 
on a curve}}
\date{1 September, 2007}
\author{Constantin Teleman\footnote
			{Partially supported by EPSRC grant GR/S06165/01}
	\and Christopher T.~Woodward\footnote
			{Partially supported by NSF grant DMS/0093647}}
\maketitle
\centerline{\emph{Dedicated to the memory of Raoul Bott}}

\maketitle
\begin{quote}
\abstract{\noindent We prove the formulae conjectured by the first 
author for the index of $K$-theory classes over the moduli stack of 
algebraic $\GC$-bundles on a smooth projective curve. The formulae 
generalise E.~Verlinde's for line bundles and have Witten's integrals 
over the moduli space of stable bundles as their large level limits.  
As an application, we prove the Newstead-Ramanan conjecture on the 
vanishing of high Chern classes of certain moduli spaces of semi-stable 
$\GC$-bundles.}
\end{quote}

\section*{Introduction}

Let $\GC$ be a reductive, connected complex Lie group and $\frM$ the 
moduli stack of algebraic $\GC$-bundles over a smooth projective curve 
$\Sigma$ of genus $g$. In this paper, we determine the \textit{analytic 
index} on a dense subring of the topological $K$-theory of $\frM$.  For 
line bundles, we recover the famous formula due to E.~Verlinde \cite
{ver}, which we extend to include the \textit{Atiyah-Bott classes}, 
described in \S1. From this angle, our index is analogous to Witten's
cohomological integration formula \cite{wit} over the moduli space of
semi-stable bundles, which appears for us in the large level limit of
the index.  Like the Verlinde formula, but unlike Witten's, our index
is expressed as a finite sum; this removes the convergence problems
and consequent regularisation in \cite{wit}. While other
regularisations have been considered in the literature \cite{jk}, ours
is intrinsically meaningful in topological $K$-theory, and expresses
the fact that indexes of vector bundles over $\frM$, and not just
those of line bundles, are controlled by finite-dimensional Frobenius
algebras \cite{tel3}.

For a smooth projective variety $X$, the analytic index of a
holomorphic vector bundle $V$ can be defined as the alternating sum
$\chi(X;V)$ of its sheaf cohomologies. This agrees with the \textit
{topological index} of $V$, defined by the Gysin map to a point in
topological $K$-theory. The construction extends to certain well-behaved 
Artin stacks. Thus, when $G$ acts on $X$, a vector bundle
over the quotient stack $X/\GC$ corresponds to an equivariant bundle
$V$ over $X$, and the holomorphic Euler characteristic $\chi(X/\GC;
V)$ (defined, say, using the simplicial bar construction of $X/G$)
agrees with the invariant part of the virtual $\GC$-representation
$\chi(X;V)$.  Now, $\chi(X;V)$ agrees with the \textit{equivariant 
topological index}, the image of $V$ under the Gysin map from $K^0_
{\GR}(X)$ to $R_G$, the representation ring of the maximal compact 
subgroup $\GR\subset\GC$. Regarding the map $R_G\to\bZ$ which extracts 
the invariant part of a representation as the Gysin projection from 
the classifying stack $B\GR$ to a point gives us an ``analytic = 
topological" index theorem for $X/G$.

Our stack $\frM$ fails a basic test for good behaviour: it has 
infinite type. Thus, when $\GC$ is a torus, $\frM$ has infinitely many 
connected components, labelled by $H^2(\Sigma, \pi_1 G)$. 
Nonetheless, $\frM$ has a distinguished \textit{Shatz stratification}. 
For a torus, the strata are the connected components. Finite, open 
unions of strata can be presented as quotients 
of smooth quasi-projective varieties by reductive groups; this allows
us to use familiar techniques of sheaf cohomology. In addition, 
special geometric features of the stratification --- 
reflected in the properties of canonical parabolic reductions of 
$G$-bundles --- ensure the finiteness of sheaf cohomology, and allows 
us to define the index, for a sub-ring of \textit{admissible} $K$-theory 
classes. When $G$ is simply connected, these classes are dense 
in the rational $K$-theory of $\frM$, in the topology induced by the 
stratification. The index is not continuous in this topology and does 
not extend to all of $K^0$; nevertheless, interesting limits do exist, 
such as in our application to the Newstead conjecture in \S\ref{8}.

This extension of Verlinde's formula, capturing the index of
vector bundles, emerged from the discovery that a certain twisted
$K$-theory was the topological home for the index of line bundles over
$\frM$ \cite[\S8]{fht}. Thus motivated, formulae for the index of
admissible classes were proposed in \cite{tel3}, equating the analytic
index over $\frM$, defined from coherent sheaf cohomology, with a
topological index defined in twisted $K$-theory. (See also the
informal notes \cite{tel4}.) The topological index can be calculated
by the Atiyah-Bott fixed-point method, and when $\pi_1G$ is free, one
obtains a formula in terms of the maximal torus $T$ and the Weyl
denominator.  In this situation, \cite{tel3} offers two conjectural
formulae for the analytic index over $\frM$ which do not involve
twisted $K$-theory: a localisation formula, which reduces the index to
the stack of principal $T$-bundles, and a Verlinde-like formula
involving a finite sum over conjugacy classes.

We prove these formulae here. It is clear that our method leads to 
the equality of analytic and topological indexes for all compact groups, 
but for simplicity we confine ourselves to \textit{connected groups with 
free $\pi_1$}; explicit formulae for more general groups require 
additional calculations (Remark~\ref{discon}).

In principle, we also solve the index problem over the more traditional 
moduli \textit{space} $M$ of semi-stable bundles. For large levels (twists 
by large line bundles), the contribution of unstable strata of $\frM$ 
vanishes, and the index over $\frM$ is equal to that over the open 
sub-stack $\frM^{\on{ss}}$ of semi-stable bundles \eqref{indexlemma}. 
The cohomology of a coherent sheaf over $\frM^{\on{ss}}$ agrees with that 
of its direct image to $M$. The index over $M$, which is a projective 
variety, depends (quasi-)polynomially on the level, so we can give a 
formula (unpleasant, but explicit) for the index of (the direct image 
to $M$ of) admissible classes, at any level. When all semi-stable 
$\GC$-bundles are stable,\footnote{When $G\neq \GL(n)$, this condition 
can only hold if we include parabolic structures.} $\frM^{\on{ss}}$ is an 
orbifold with coarse quotient $M$, and the rational cohomology 
calculation of \cite{ab} shows that we generate all of $K^0(M;\bQ)$ 
in this way.  

Our proof relies on a remarkable symmetry of the index over $\frM$ 
which is \emph{absent} on $M$, or on any finite-type approximation. 
The symmetry arises from a loop group version of  Bott's reflection 
argument \cite{bott}, a \textit{Hecke correspondence}. (This device 
was already used in \cite{bs} in relation to the Verlinde formula.) 
The reader should refer to \S\ref{atbott} below for the definitions in
what follows.  Choose an \textit{admissible line bundle} $\cL$ and an
\textit{index bundle} $E^*_\Sigma V$. For a weight $\lambda$ of the
maximal torus $T$, denote by $V_\lambda$ the holomorphically induced
virtual $\GC$-representation.  Regard the index of
$\cL\otimes\exp[tE^*_\Sigma V]\otimes E_x^*V_ \lambda$, a formal
series in $t$, as the $e^\lambda$-coefficient for a Fourier series on
$T$, with values in $\bQ[[t]]$. This series turns out to be
anti-invariant for a certain action of the affine Weyl group, and is
thereby constrained to represent a sum of $\delta$-functions at
prescribed, \emph{regular} points of $T$.  Regularity of its support,
combined with Atiyah's localisation theorem for the index of
transversally elliptic operators, implies that the index distribution
only sees the contribution of principal bundles whose structure group
reduces to $T$. That can be calculated by Riemann-Roch, leading to our
explicit index formula.
 
The paper is organised as follows. In \S\ref{2} we describe the admissible 
$K$-classes and define their analytic index. We include a brief review 
of the stratification of $\frM$ and the local cohomology vanishing 
results of \cite{tel2}. Section \ref{3} contains the precise statements 
of our formulae. The proof is split into \S\ref{abred}, where we check 
the anti-symmetry of the index distribution, and \S\ref{5} where we 
eliminate the contributions of non-toric principal bundles. 

The last sections contain two applications. In \S\ref{witten}, we show
how Witten's integration formulae over $M$ arise from our index
formula in the large level limit; we only give full details for
$\SL(2)$. (The formulae were proven for $\SL(r)$ by Jeffrey-Kirwan
\cite{jk} and, independently of our work but simultaneously, by
Meinrenken \cite{mein} for compact, $1$-connected $\GC$.) Section~\ref
{kaehler} enhances our index formulae by incorporating K\"ahler
differentials, needed in our next application in \S\ref{8} to a
conjecture of Newstead and Ramanan. The original version, proved by
Gieseker \cite{gies}, asserted the vanishing of the top $2g-1$ Chern
classes of the moduli space of stable, odd degree vector bundles of
rank $2$ on $\Sigma$. An analogue in rank $3$ was settled by Kiem and
Li \cite{kiemli}. We generalise this to the vanishing of the top
$(g-1)\ell$ rational Chern classes of the moduli space $M$ of stable
principal bundles with semi-simple structure group of rank $\ell$,
whenever $M$ (or a variant decorated with parabolic structures) is a
compact orbifold.

The appendix reviews some background on the topological $K$-theory
of $\frM$ and on its variants decorated with parabolic structures; 
the exotic parabolic structure associated to the simple affine 
root of $\frg$ plays a special r\^ole in the proof. We do not 
review general properties of stacks and their cohomology, these 
matters having had increasing coverage in the literature since 
the detailed treatments \cite{bl, ls}; a review suited to our 
needs is found in \cite{tel1, tel2}.

We thank the referee for a careful reading of the manuscript and for
many helpful suggestions.

\subsection*{Notation.} $\GC$ is a reductive group, $T$ a maximal 
torus and $B\supset T$ a Borel subgroup. $\GR, T_k$ will be the compact 
forms and $\frgR, \frt_k$ their Lie algebras. The co-weight lattice 
of $T$, $\log(1)/2\pi\mi$, lies in $\mi\frt_k$; its $\bZ$-dual is the weight
lattice in $\mi\frt_k^\vee$. $W$ is the Weyl group and $\Delta:= 
\prod_{\alpha>0}2 \sin (\mi\alpha/2)$ the Weyl denominator. The 
Weyl vector $\rho$ is the half-sum of the positive roots. The simple 
roots are $\alpha_1,\dots,\alpha_\ell$; when $\frg$ is simple, the 
simple affine root $\alpha_0$ sends $\xi\in\frt$ to $1-\hroot(\xi)$, 
with the highest root $\hroot$. (The affine root vector of $\alpha_0$ 
is $z^{-1}e_\hroot$.) The representation ring of $\GR$ is denoted by 
$R_G$, and $\bC{R}_G := \bC\otimes R_G$.

\section{Atiyah-Bott classes}
\label{2}

In this section, we introduce the Atiyah-Bott classes and admissible 
classes. We then define their analytic index and derive its finiteness 
from the local cohomology vanishing results of \cite{tel2}. This requires 
a brief review of the Shatz stratification.

\subsection{Admissible classes.} \label{atbott}
Given a representation $V$ of $\GC$, call $E^*V$ the vector bundle 
over $\Sigma\times\frM$ associated to the universal $\GC$-bundle. 
Call $\pi$ the projection along $\Sigma$, $\sqrt K$ a square root 
of the relative canonical bundle, and $[C]$ the topological $K_1
$-homology class of a $1$-cycle $C$ on $\Sigma$. Consider the 
following  classes in the topological $K$-theory of $\frM$:

\begin{enumerate}
\item[(i)] The restriction $E^*_xV \in K^0(\frM)$ of $E^*V$ to a point
$x\in\Sigma$;
\item[(ii)] The slant product $E^*_CV:=E^*V/[C]\in K^{-1}(\frM)$ of
$E^*V$ with $[C]$;
\item[(iii)] The Dirac index bundle $E^*_\Sigma V:=R\pi_*
(E^*V\otimes\sqrt K)\in K^0(\frM)$ of $E^*V$ along $\Sigma$;
\item[(iv)]  The inverse \emph{determinant of cohomology}, 
$D_\Sigma V:= \det^{-1} E^*_\Sigma V$. 

\end{enumerate}
We call the classes (i)--(iii) the \textit{Atiyah-Bott generators}; 
they are introduced in \cite[\S2]{ab}, along with their counterparts 
in cohomology, and can also be described from the K\"unneth decomposition 
of $E^*V$ in 
\[
K^0(\Sigma\times \frM) \cong K^0(\Sigma)\otimes K^0(\frM) 
			\oplus K^1(\Sigma)\otimes K^1(\frM),
\]
by contraction with the various classes in $\Sigma$.  
Classes (i) and (iv) are represented by algebraic vector 
bundles, while (iii) can be realised as a perfect complex of 
$\cO$-modules. The class $E^*_CV$ in (ii) is not algebraic. 
Note that $\det E^*_\Sigma V =  \det {R}\pi_*(E^*V)$ when  
$\det{V}$ is trivial; an important example is the canonical bundle 
$\cK = \det E^*_\Sigma\frg$ of $\frM$, defined from the adjoint 
representation $\frg$.

For general (non-simply connected) groups, determinant line bundles 
are quite restrictive; we will consider more generally line bundles 
which have a \emph{level}, defined below, and call them \emph
{admissible} if their level exceeds that of $\cK^{1/2}$. (Recall 
\cite{ls} that $\cK$ has a distinguished \textit{Pfaffian} square 
root.) Products of an admissible line bundle and any number of 
Atiyah-Bott generators span the ring of \textit{admissible classes}. 

\subsection{Line bundles with a level.} \label{level}
To certain line bundles on $\frM$ we will now associate a {\em level}, 
a quadratic form on the Lie algebra $\frg$. Briefly, for any $V$, 
the level of $\det E^*_\Sigma{V}$ is the trace from $\xi,\eta\mapsto 
\Tr_V(\xi\eta)$, and we wish to extend this by linearity in the first 
Chern class of the line bundle.

Riemann-Roch along $\Sigma$ expresses $c_1({E}^*_\Sigma{V})$ as the 
image of $ch_2(V) = \frac{1}{2} c_1^2(V) - c_2(V)$ under \emph{
transgression along $\Sigma$}, $\tau: H^4(BG;\bQ) \to H^2(\frM;\bQ)$ 
(construction (\ref{atbott}.iii) in cohomology). It is important 
that $\tau$ is injective (Remark \ref{injtrans}). We now identify 
$H^4(BG;\bR)$ with the space of invariant symmetric bilinear forms 
on $\frgR$ so that $\Tr_V$ corresponds\footnote{It is more standard 
to identify $\Tr$ with $2\,ch_2$; our choice here avoids factors 
of $2$ elsewhere.} to $ch_2(V)$. We say that the line bundle $\cL$ 
\emph{has a level} if its Chern class $c_1(\cL)$ agrees with some 
$\tau(h)$ in $H^2(\frM; \bQ)$; the form $h$, called \emph{level} 
of $\cL$, is then unique. 

For $\SL_n$, the level of the positive generator of $\mathrm{Pic}(\frM)$ 
is $-\Tr_{\bC^n}$ in the standard representation; the calculation 
is due to Quillen. For another example, the level of $\cK^{-1/2}$ is 
$c:= -\frac{1}{2}\Tr_\frg$. \emph{Positivity} of a level refers to 
the quadratic form on $\frgR$; thus, $D_\Sigma V$ has positive 
level iff $V$ is $\frg$-faithful. Finally, $\cL$, with level $h$, 
is \emph{admissible} iff  $h>-c$ as a quadratic form.

\begin{remark}
\begin{trivlist}\itemsep0ex
\item (i) When $G$ is simply connected, $\tau: H^4(BG;\bZ) \to 
H^2(\frM;\bZ)$ is an isomorphism, but this fails (even rationally) 
as soon as $\pi_1{G} \neq 0$. Line bundles with a level satisfy a 
prescribed relation between their Chern classes over the different 
components of $\frM$: cf.~\eqref{first}.
\item (ii) The traces span the negative semi-definite cone in $H^4(BG;\bR)$; 
so $\cL$ has positive level iff $c_1(\cL)$ lies in the $\bQ_+$-span 
of the $c_1(D_\Sigma V)$'s, with $\frg$-faithful $V$.
\item (iii) For semi-simple $G$, $\cK$ has negative level, and so $\cO$ 
is admissible. This fails for a torus, but positive-level line bundles 
are admissible for any $G$.
\item (iv) For $g>1$ and simply connected $G$, positivity of the 
level is equivalent to amplitude on the moduli space. (Suffices to check 
this for simple $G$: recall then that $\mathrm{Pic}(\frM)=\bZ$ and that 
$\cK^{-1}$ is ample.) When $\pi_1G \neq 0$, the positive level 
condition is much more restrictive.  
\end{trivlist}\end{remark}

\subsection{The index of an admissible class.}
We first recall the finiteness result which enables us to define the
index of admissible classes by means of sheaf cohomology. It is a 
consequence of \cite {tel2}, combining the relative case of the main 
theorem in \textit {loc.\ cit.}, \S5 with the discussion of $\frM$ 
in \S8 and \S9. For the reader's convenience we will also outline 
the proof in \S\ref{vanish}, after we review the stratification 
of $\frM$.

Let $\cE$ be the twist of an external tensor product $\boxtimes\, E^*V_k$ 
of universal bundles  over $\Sigma^n\times \frM$ by an admissible line 
bundle $\cL$. Call $\underline\cE$ the direct image to $\Sigma^n\times 
M$, the moduli space, of the restriction of $\cE$ to the semi-stable part 
$\frM^{\on{ss}}$. Consider the projections $\phi$ and $\underline\phi$ from 
$\Sigma^n\times\frM$, resp.\ $\Sigma^n\times{M}$, to $\Sigma^n$. 

\begin{lemma}\label{indexlemma}
The total direct image $\bigoplus_i R^i\phi_*\cE$ is coherent  on 
$\Sigma^n$. For large enough $\cL$, depending on the $V_k$, it 
agrees with $\bigoplus_i R^i\underline{\phi}_*\underline\cE$. \qed
\end{lemma}
\noindent A lower bound for the level of $\cL$ can be given, linear in 
the highest weights of the $V_k$, see \S\ref{vanish}. 

Choose now cycles $C_k$ on $\Sigma$ of various dimensions, but with 
even total degree. We wish to define the index of $\cL\otimes 
\bigotimes_k E^*_{C_k}V_k$ over $\frM$ as the Euler characteristic 
of its coherent sheaf cohomology, but the non-algebraic classes 
(\ref{atbott}.ii) impose the indirect   

\begin{definition}\label{indexdef}
The \textit{index} over $\frM$ of $\cL\otimes\bigotimes_k E^*_{C_k}V_k$ 
is the pairing of $\times_k [C_k]\in K_0(\Sigma^n)$ with the 
topological $K$-theory class of $\sum_i(-1)^iR^i\phi_*\cE$. 
\end{definition}
\noindent When all the $C_k$ are even, we can switch and push down 
along $\Sigma^n$ first, recovering the Euler characteristic. 
When $G$ is abelian, the index theorem applied to the components 
of $\frM$ shows that our index only depends on the underlying 
topological $K$-class of the bundle. That is not so obvious in general, 
but will follow for instance from our abelianisation formula 
\eqref{local}.

\subsection{Shatz stratification.} \label{shatz}
Recall that any $\GC$-bundle over $\Sigma$ admits a canonical
reduction of structure group to a standard parabolic subgroup $P$ of
$G$, for which the associated bundle with Levi structure group is
semi-stable.  Topologically, this reduction is classified by a co-weight 
of $P/[P,P]$; we identify this with a (possibly fractional) dominant
co-weight $\xi$ of $\frg$, called the \textit{instability type} of the
original bundle. Then, $P$ is the standard parabolic subgroup defined
by $\xi$; we will denote it by $P_\xi$ and its Levi subgroup by
$G_\xi$. If $\frM_\xi$ denotes the stack of $G$-bundles of type $\xi$,
we have an algebraic stratification \cite{sh:de, ab}
\[
\frM = \bigcup\nolimits_\xi \frM_\xi.
\]
Sending a $P_\xi$-bundle to its associated Levi bundle gives a
morphism from $\frM_\xi$ to the stack $\frM_{\GC_\xi,\xi}^{\on{ss}}$ of
semi-stable principal $\GC_\xi$-bundles of type $\xi$; the fibres are
quotient stacks of affine spaces by nilpotent groups. The virtual 
normal bundle for the morphism $\frM_{\GC_\xi,\xi}^{\on{ss}} \to\frM$ 
is the complex 
\[
\nu_\xi = R\pi_* E^*(\frg/\frg_\xi)[1].
\]
Its $K$-theory Euler class should be the alternating sum of exterior 
powers\footnote{Recall that the $p$th exterior power $\lambda^p$ of a complex 
$V^0  \xrightarrow{\delta} V^1$ is the complex with $q$th space $\Lambda^
{p-q}V^0\otimes \mathrm{Sym}^qV^1$ and obvious differential induced 
by $\delta$. A similar definition applies to symmetric powers.} 
\[
\lambda_{-1}(\nu_\xi^\vee):= \sum (-1)^p\lambda^p(\nu^\vee_\xi),
\] 
but for now this infinite sum is only a formal expression, whose meaning 
is to be spelt out. 

\subsection{Local cohomology.}\label{spectral}
Finite, open unions of Shatz strata
\[
\frM_{\leq \xi} = \bigcup\nolimits_{\mu \leq \xi} \frM_\mu
\]
can be presented as quotient stacks of smooth quasi-projective varieties 
by reductive groups. The cohomology with supports over $\frM_\xi$ of a 
vector bundle $\cE$ is
\begin{equation}\label{localcoh}
H^\bullet_{\frM_\xi}(\frM_{\leq \xi},\cE_{\leq \xi}) =
H^{\bullet+d_\xi}(\frM_\xi, \cR_\xi \cE) 
\end{equation}
where $d_\xi$ is the co-dimension of $\frM_\xi$ and $\cR_\xi \cE \to 
\frM_\xi$ the sheaf of ``$\cE$-valued residues along $\frM_\xi$", 
the cohomology sheaves relative to the complement of $\frM_\xi$. 
Pushing down to $\frM_{\GC_\xi,\xi}^{\on{ss}}$ and passing to the 
associated graded sheaf for the filtration by order of the pole 
leads to
\begin{equation}\label{gradedcoh}
 {H}^\bullet(\frM_{\xi,\GC_\xi}^{\on{ss}}, 
	\cE_\xi \otimes\on{Eul}(\nu_\xi)^{-1}_+)
\end{equation}
where $\cE_\xi$ is the restriction to $\frM_{\GC_\xi,\xi}^{\on{ss}}$, 
while the complex of sheaves
\[
\on{Eul}(\nu_\xi)^{-1}_+ := \Sym R\pi_*
	\left(E^*(\frp_\xi/\frg_\xi)[1]^\vee \oplus 
	R\pi_* E^* (\frg/\frp_\xi)[1]\right)
	\otimes \det(R\pi_* E^* (\frg/\frp_\xi)[1]) [d_\xi]
\]
is formally an inverse to the Euler class $\lambda_{-1}(\nu_\xi^\vee)$ 
which ``prefers" the $\xi$-negative eigenvalues in the geometric 
expansion.  

\subsection{Finiteness and vanishing.}\label{vanish}
All eigenvalues of $\xi$ appearing in $\on{Eul}(\nu_\xi)^{-1}_+$ are 
negative, with finite multiplicity. The determinant factor has weight 
$c(\xi,\xi)$ (negative, as $\xi\in\mi\frt_k$). An admissible line 
bundle factor $\cL$ in $\cE$ changes this behaviour to $(h+c)(\xi,\xi)$. 
Atiyah-Bott bundles $E^*_?V$ alter this behaviour linearly in $\xi$. 
Overall, for any admissible $\cE$, the $\xi$-invariant part of 
$\cE\otimes \on{Eul}(\nu_\xi)^{-1}_+$ is finite-dimensional, and 
vanishes for all but finitely many $\xi$.

It follows that almost all cohomologies \eqref{localcoh} vanish, and 
the index of $\cE$ over $\frM$ is the sum of finitely many local 
contributions over the $\frM_\xi$. Passage to the $\on{Gr}$ does not 
change the index and we obtain  
\begin{equation} \label{notlocalisation}
\Ind(\frM,\cE) = \sum\nolimits_{\xi} \Ind(\frM_{\GC_\xi,\xi}^{\on{ss}}, 
		\cE_\xi\otimes \on{Eul}(\nu_\xi)^{-1}_+).
\end{equation}
Lemma~\ref{indexlemma} is the relative version of this story for the 
projection $\phi$ to $\Sigma^n$, with $R^i\phi_*$ replacing cohomology 
and $\sum (-1)^i R^i\phi_*\in {K}^0(\Sigma^n)$ replacing the index.

\begin{remark}  Formula \eqref{notlocalisation} is related to the 
\textit{non-abelian localisation} principle of Witten \cite{wit}.
When presenting $\frM_{\leq \xi}$ as a quotient of a manifold by 
a reductive group, the $\bar\partial$ operator can be deformed 
so that the invariant part of its kernel localises at the critical 
points of the norm-square of the moment map, leading to the individual 
contributions in \eqref{notlocalisation}, see Paradan \cite{pa:lo}.
\end{remark}

\begin{remark} Inadmissible $\cL$'s can have infinitely many 
contributions to \eqref{notlocalisation}.  However, 
when $\GC$ is semi-simple, the theorems of Kumar \cite{ku:km} and 
Mathieu \cite{ma:km} imply the vanishing of all the direct images 
for negative $\cL$: $\frM$ is then isomorphic to a quotient 
of the generalised flag variety $X:= \GC\pz/\GC[[z]]$ for the loop 
group by the group $G[\Sigma \setminus\{ x \}]$ of algebraic maps 
on the punctured curve, and the cohomology of $\cL$ vanishes over $X$.
\end{remark}

\section{The index formulae}
\label{3}

The index formulae involve a sum over deformations of certain 
\emph{Verlinde conjugacy classes} in $G$, which appear in the 
formula for line bundles. We start by recalling that story.  

\subsection{Isogenies from admissible levels.}\label{adm}
Contraction $\xi\mapsto \iota(\xi)h$
with the level $h$ of an admissible line bundle $\cL$ maps the
co-weight lattice to its dual, the weight lattice. This map descends
to a homomorphism $\chi: T \to T^\vee$, the dual torus. The
homomorphism $\chi'$ defined from $h':=h+c$ is an isogeny, with kernel
$F\subset T$. Let $F_\rho$ be the translate of $F$ lying over
$e^{2\pi\mi\rho}\in T^\vee$.  This last point does not depend on the
Weyl chamber used to define $\rho$, and gives the Spin covering of $T$
in the adjoint representation $\frg$.

\begin{example}  If $G = \SL(n)$, $\mathrm{Pic}(\frM)\cong\bZ$,
with positive generator $\cO(1)= D_\Sigma \bC^n$ and $\cK=\cO(-2n)$. 
$T^\vee$ is the maximal torus of $\bP\SL(n)$ and $\chi$, for $\cO(1)$, 
is the natural projection.  Hence, for $\cL=\cO(l)$, $F=F_\rho$ 
comprises the $(l+n)^{th}$ roots of the centre of $\SL(n)$. The 
analogue holds for simply connected, simply laced groups, if $n$ 
is replaced by the dual Coxeter number.
\end{example}

A formula of E. Verlinde (first given in the context of conformal
field theory) describes the index of a determinant line bundle over
$\frM$.  Let $\Theta$ be the sum of delta-functions
on the \emph{regular} $\GR$-conjugacy classes through $F_\rho$,
divided by the order $|F|$ of $F$. Define a linear map $R_G\to\bZ$ on
representations by
\begin{equation}\label{orig}
U \mapsto \Theta(U) = \int_{\GR}{\Tr_U(g)\cdot\Theta(g)dg} = 
	\sum_{f\in F_\rho^\mathrm{reg}/W}\Tr_U(f)\cdot\frac{\Delta(f)^2}{|F|}.
\end{equation}
(Recall that we normalised the Weyl denominator so that $\Delta(f)^2$ 
is the volume of the conjugacy class.)  Let $\theta(f) = \Delta(f)^2/|F|$; 
Verlinde's formula\footnote{See, e.g.~\cite{amw} for semi-simple $\GC$; 
we shall reprove it below when $\pi_1G$ is free.} gives the index of $\cL$ 
as
\begin{equation}\label{oldverlinde}
\Ind(\frM;\cL) = 
	 \sum_{f\in F_\rho^\mathrm{reg}/W} \theta(f)^{1-g}.
\end{equation}

\begin{remark}\label{spinsign}
There is a version of formula \eqref{oldverlinde} with $F_\rho$ replaced 
by $F$. The components $\frM^{(\gamma)}$  of $\frM$ are labelled by $\gamma\in
\pi_1\GC$, and the Spin covering of the adjoint representation of $\GC$ 
defines a character $\sigma:\pi_1G \to \{\pm1\}$. The calculations of 
\S\ref{5} give a \emph{graded} index formula
\begin{equation}\label{gradedindex}
\sum_{\gamma\in\pi_1G} \sigma(\gamma)\cdot \Ind(\frM^{(\gamma)};\cL) = 
	 \sum_{f\in F^\mathrm{reg}/W} \theta(f)^{1-g}.
\end{equation}
The same applies to our generalised index formulae below. 
\end{remark}

\begin{remark}
The kernel of the pairing $(U,U') \mapsto \Theta(U\otimes U')$ is the 
ideal $I_h\subset R_G$ of virtual characters which vanish on $F_\rho^
\mathrm{reg}$. We obtain a non-degenerate pairing on the quotient $R_G/I_h$, 
which becomes an integral Frobenius algebra, the \emph {Verlinde ring 
at level~$h$}. Its complex spectrum is $F_\rho^\mathrm{reg}/W$. A folk 
result asserts that a Frobenius algebra is the same as a $2$-dimensional 
topological field theory, and formula \eqref{oldverlinde} is the 
``partition function" for a genus $g$ surface in the Verlinde 
ring.
\end{remark}

\subsection{Deformations.}
Given a representation $V$ of $\GC$, consider the following formal
one-parameter family of transformations on $\GC$:
\begin{equation}\label{coord}
	g\mapsto m_t(g):= g\cdot\exp\left[t\,\nabla\Tr_V(g)\right],
\end{equation}
with the gradient in the bilinear form $h'$. This descends to the space 
$\GC/\mathrm{Ad}\GC$ of conjugacy classes; note from the Ad-invariance 
of $\Tr_V$ that singular classes remain singular. Restricting to 
conjugacy classes in $\GR$ and composing with $\Theta$ gives a formal 
$t$-family $\Theta_t:= \Theta\circ m_t$ of distributions, even though the 
points $f_t$ of its support, the solutions to $m_t(f_t) = f$, can move 
in complex directions (this happens when $\Tr_V$ is not real):
\[
\Theta_t(U) := \int_{\GR}{\Tr_U(g)\cdot\Theta(m_t(g))\,dg}
	= \sum_{f\in F_\rho^\mathrm{reg}/W}\Tr_U(f_t)
	\cdot\theta_t(f_t)\in\bC[[t].
\] 
The $\theta_t$ are described as follows. Call $\hess_V$ the Hessian 
of $\Tr_V$: $\hess_V(u)(\xi,\eta)= \Tr_V(u\xi\eta)$, for $u\in T$ and 
$\xi,\eta \in\frt$, and denote by $\hess_V(u)^\dagger$ its conversion via 
$h'$ to an endomorphism of $\frt$. In view of the volume scaling under 
\eqref{coord}, we have 
\begin{equation}\label{thetas}
\theta_t(f_t) = \det\nolimits^{-1}\left[1+t\hess_V(f_t)^\dagger\right]
		\cdot\frac{\Delta(f_t)^2}{|F|}.
\end{equation}
Here is an alternate description of $\theta_t$. Strictly speaking, it 
applies only when $\Tr_V$ is real; see the closely related Fourier 
expansion of $\Theta_t$ in $\S\ref{sum}$, which is free of this flaw.  
The push-down of $\Theta$ to the space $G_k/\mathrm{Ad}G_k = T_k/W$ 
of unitary conjugacy classes is
\begin{equation} \label{torusdelta}
\sum\nolimits_{F_\rho/W} \theta(f)\cdot\delta_f = 
	\Delta^2 \cdot\delta_\rho\circ\chi',
\end{equation}
with the delta-function $\delta_\rho$ at $e^{2\pi\mi\rho}\in T^\vee$
and the isogeny $\chi': T\to T^\vee$ of \S\ref{adm}. Viewing  $d\,\Tr_V$ 
as a map $T \to\frt^\vee$, $\chi'$ has a formal deformation to 
$\chi'_t := \chi'\cdot\exp(t d\,\Tr_V)$, under which \eqref{torusdelta} 
deforms to
\begin{equation}\label{deformtheta}
\sum \theta_t(f_t)\cdot\delta_{f_t} = \Delta^2\cdot\delta_\rho\circ \chi'_t.
\end{equation}

\begin{remark}\label{twisting}
Pulling back Fourier modes on $T^\vee$ by $\chi'_t$ defines a group 
homomorphism from $\pi_1T$ to the units in $\bC{R}_T[[t]]$. This defines a 
\emph{(higher) twisting} for the equivariant $K$-theory $K_T\left(T; 
\bC[[t]]\right)$. This extends to a twisting for the conjugation action 
of $G$ on itself, and the twisted $K$-group $K_G(G)$ turns out to be 
the quotient of $\bC R_G[[t]]$ by the kernel of the pairing $\Theta_t
(U\otimes U')$. It is a Frobenius algebra over $\bC[[t]]$, deforming 
the complex Verlinde ring at $t=0$. See \cite{tel3} for more details. 
\end{remark}

\subsection{Even Atiyah-Bott generators.}
We incorporate the index bundles (\ref{atbott}.iii) into Verlinde's 
formula by means of a generating function
\[
\cL\otimes\exp[t_1E^*_\Sigma V_1+\ldots +t_nE^*_\Sigma V_n]\otimes
			E^*_xU \in K^\bullet(\frM)[[t_1,\ldots,t_n]].
\]
Let $\theta_\mathbf{t}(f)$ denote the multi-parameter version of 
$\theta_t(f)$, for $\mathbf{t} = (t_1,\ldots,t_n)$.

\begin{theorem}[Index formula for even classes]\label{main}
\[
\Ind\left(\frM;\cL\otimes\exp[t_1E^*_\Sigma V_1+\ldots
			+t_nE^*_\Sigma V_n]\otimes E^*_xU\right) = 
			\sum_{f\in F_\rho^\mathrm{reg}/W}
			{\theta_\mathbf{t}(f)^{1-g}\cdot 
					\Tr_U(f_\mathbf{t})}.
\]
\end{theorem}
\noindent With $\mathbf{t}=0$ and the trivial representation $U$, this
recovers \eqref{oldverlinde}.

\subsection{Example.} \label{su2}
When $G=\SL(2)$, $\cL = \cO(l)$ and $\Tr_V =\sum\varphi_nu^n$ on matrices 
with eigenvalues $\{u,u^{-1}\}$, we have, as conjectured 
in \cite{tel3},
\[
\Ind\left(\frM; \cL\otimes \exp[tE^*_\Sigma V]\right) = 
		\sum\limits_{\zeta_t} {\left[\frac{2l + 4 + t\ddot \varphi
		(\zeta_t)}{\left|\zeta_t-\zeta_t^{-1}\right|^2}\right]^{g-1}}
\]
where the $\zeta_t$ range over the solutions of $\zeta_t^{2l+4}\cdot\exp
\left(t\dot\varphi(\zeta_t)\right) = 1$ with positive imaginary part, $\dot 
\varphi(u) = \sum n\varphi_nu^n$ and $\ddot\varphi(u)=\sum n^2\varphi_nu^n$.

\subsection{Odd generators.}  \label{odd}
The bilinear form $h'+ t\hess_V(u)$ on $\frt$ is non-degenerate;  
denote by $\langle\: |\: \rangle(u)$ the inverse form on $\frt^\vee$. 
To an even product $\psi$ of odd Atiyah-Bott generators (\ref{atbott}.ii), 
we assign a function $[\psi](u)$ on $T$ as follows: split $\psi$ into 
quadratic factors $E_C^*U\wedge E_{C'}^*U'$, replace each factor by the number
\[
-\#(C\cap C')\cdot\langle d\,\Tr_U(u)| d\,\Tr_{U'}(u)\rangle(u),
\]
where $\#(C\cap C')$ is the intersection pairing, and sum over all 
possible quadratic splittings, with signs as required by re-ordering.  
Set $[\psi](u)=0$ if $\psi$ is odd. We shall see in \S\ref{5} that 
$[\psi](u)$ is expressible in terms of the integral of the Chern 
character $\on{Ch}(\psi)$ against a Gaussian form $\exp\{[h'+tH_V(u)]
\otimes\eta\}$ on the Jacobian of $T$-bundles on $\Sigma$; in particular, 
it only depends on the $K$-theory class of $\psi$. The following gives 
the index for odd and even classes; for simplicity we use a single $V$.

\begin{theorem}[Index formula for general classes] \label{general}
\[
\mathrm{Ind}\left(\frM;\cL\otimes \exp[tE^*_\Sigma V]\otimes 
		E^*_xU\otimes\psi\right)= 
		\sum_{f\in F_\rho^\mathrm{reg}/W} \Tr_U(f_t)\theta_t(f_t)^{1-g}
		\cdot[\psi](f_t).
\]\end{theorem}

\subsection{Abelianisation.}\label{abelsubsect}
We will derive our index formulae from a more conceptual ``virtual 
localisation" to the stack $\frM_T$ of $T$-bundles. Let $\nu:= 
R\pi_*E^*(\frg/\frt)[1]$ be the virtual normal bundle for the 
morphism $j:\frM_T \to \frM$. In \S\ref{eulerinverse}, we will 
see that the $K$-theoretic Euler class  $\lambda_{-1}(\nu^\vee)$ is 
well-defined after inverting the Weyl denominator and equals 
\[
\lambda_{-1}(\nu^\vee) = (-1)^{2\rho(\gamma)}\Delta^{2g-2}\cK^{1/2}
\]
on the component of $\frM_T$ of topological type $\gamma$.

\begin{theorem}\label{local} 
For admissible $\cE$, 
$\Ind\left(\frM;\cE\right) = |W|^{-1} ``\Ind"\left(
		\frM_T;j^*\cE\otimes \lambda_{-1}(\nu^\vee)^{-1}\right).$
\end{theorem}
\comment{check duals}
\noindent The right-hand side needs clarification. Each component of 
$\frM_T$ is the product of the classifying stack $BT$ with a Jacobian 
of $T$-bundles, and the index over $\frM_T$ should be the sum of the 
$T$-invariant parts of indexes over these Jacobians. Because of the Weyl 
denominator, $\lambda_{-1}(\nu^\vee)$ is not invertible in
$K^\bullet(\frM_T)$; the index of $j^*\cE\otimes \lambda_{-1}(\nu^\vee)
^{-1}$ over each Jacobian lands in $R_T[\Delta^{-1}]$, and its 
$T$-invariant part is not \textit{a priori} well-defined.  However, 
we will see in \S\ref{sum} that summing over all Jacobians leads to a 
well-defined distribution on the regular part of $T_k$, supported 
on $F_\rho^\mathrm{reg}$. We \emph{declare} the index over $\frM_T$ 
to be the invariant part (= integral over $T_k$) of this index 
distribution, \emph{after extension by zero to the singular locus}. 

A formula which does not require inverting $|W|$ will be given in 
Proposition~\ref{abel}. In a family of curves, the index is replaced 
with a $K$-theory class on the base, and the alternate formula loses 
slightly less torsion. We hope to return to this in future work.

\section{Affine Weyl symmetry} 
\label{abred}

In this section, we establish the anti-symmetry of the index under 
an action of the affine Weyl group (Proposition~\ref{waff}). This 
constrains the general form of the answer (Corollary~\ref{delta}).

\subsection{Affine Weyl action.} \label{waffaction} 
Define a group homomorphism from the co-root lattice $\Pi$ to the  
units in  $\bC R_T[[t]]$ by
\[
\Pi\ni\gamma \mapsto
\exp\left[\iota(\gamma)h' + t\frac{\partial\Tr_V}{\partial\gamma} \right].
\]
(This is the homomorphism  mentioned in Remark~\ref{twisting}.) 
Multiplication by these units combines with the Weyl transformations
into an action of the \textit {affine Weyl group} $\waff:= W\ltimes
\Pi$ on $\bC R_T[[t]]$. This action extends to the space of formal (unrestricted) Fourier series on $T$ with coefficients in $\bC[[t]]$.  

For a weight $\mu$ of $B$, call $V_\mu$ the \emph{holomorphically 
induced} virtual representation of $\GC$, that is, the $G$-equivariant 
index of the \textit{weight line bundle} $\cO(\mu)$ over the flag 
variety $\GC/B$. Define the following \textit{index series} on $T$, 
a formal Fourier series with coefficients in $\bC[[t]]$:

\begin{equation} \label{symmetry} 
\cI:= \sum\nolimits_\mu \Ind(\frM;\cL\otimes
 	\exp[tE^*_\Sigma V]\otimes E^*_xV_{\mu+\rho} ) e^{-\mu}.
\end{equation}

\noindent We use a single $t$ and only even classes to keep the notation 
manageable, but this restriction is not necessary (cf.~\S\ref
{oddserre} below). 

\begin{proposition}\label{waff}
The index series $\cI$ is anti-invariant under the affine Weyl
action. 
\end{proposition}

\begin{proof} 
Weyl anti-invariance being clear from the holomorphic induction step,
it suffices to confirm, for each simple factor of $\frg$, the sign
change under the following affine reflection $S$: the highest Weyl
reflection $s_0$, followed by subtraction of the co-root $H$ of the
highest root $\hroot$.  This is the affine analogue of the famous Bott
reflection \cite{bott}.

The stack $\frM(x,\cB)$ of $\GC$-bundles over $\Sigma$ with $B$-reduction 
at $x$ is a $G/B$-fibre bundle over $\frM$ (Example \ref{fibres}.i) and 
carries natural extensions of the weight line bundles $\cO(\mu)$ on 
the fibre. From the Borel-Weil-Bott and Leray theorems, the $e^\mu$ 
Fourier coefficient is
\[
\cI(e^\mu) = \Ind\left(\frM(x,\cB);\cL(\mu+\rho)\otimes
	\exp[tE^*_\Sigma V]\right),
\]
where $(\mu)$ is the twist by $\cO(\mu)$ and, abusively, $\cL$ 
stands for its own lift to $\frM(x,\cB)$.

Let $\frM'$ be the stack of $\GC$-bundles with parabolic structure 
at $x$ defined by the simple affine root $\alpha_0$ (Example \ref
{fibres}.iii). We have $\frM(x,\cB)\cong\frP\times_{\PSL(2)}\bP^1$, 
for the principal $\PSL(2)$-bundle $\frP\to \frM'$ determined by 
$\alpha_0$. Call $p$ the projection to $\frM'$. For a vector bundle 
$\cE$ over $\frM(x,\cB)$ which tensors into $\SL(2)$-equivariant 
bundles over the two factors,\footnote{But such that the diagonal 
action factors through $\PSL(2)$.} define a new bundle $\bD\cE$ by 
dualising the $\bP^1$-factor and then twisting by the relative canonical 
bundle $\cO(-\hroot)$. Relative Serre duality along $\bP^1$ gives 
$R^i p_*\bD\cE = R^{1-i}p_*\cE$ (as can be seen from $\SL(2)$-equivariant 
Serre duality on $\bP^1$ and self-duality of $\SL(2)$-representations).  
Integrating over $\frM'$ shows that the indexes of $\cE$ and $\bD\cE$ 
over $\frM(x,\cB)$ differ by a sign, and we will prove our proposition 
by relating $S$ to $\bD$. 

We claim that $\cL$ factors for the fibre product presentation of 
$\frM(x,\cB)$ as  
\begin{equation}\label{factorline}
\cL \cong \cL(-{\textstyle\frac{1}{2}}\iota(H)h)\boxtimes 
\cO({\textstyle\frac{1}{2}}\iota(H)h). 
\end{equation}
Then, $\bD\cL = \cL(-\iota(H)h - \hroot)$. Further,\footnote{Both 
sides are parallel to $\hroot$, so the equality only needs testing 
against $H/2$, when the two sides become $c(H,H)/2$ and $\rho(H) +1$, 
which are equal to the dual Coxeter number.}  $\iota(H)c = \rho - s_0\rho 
+\hroot$, and we get 
\begin{equation}\label{dualine}
	\bD\left[\cL(\mu+\rho)\right] = \cL(s_0\mu-\iota (H)h'+\rho) = 
		\cL(S\mu + \rho),
\end{equation}
confirming the proposition for $t=0$. 

To verify \eqref{factorline}, note that \textit{some} such 
formula must hold, with $h$ replaced by a fixed multiple of itself; 
namely, the one which renders the first factor trivial along the fibres 
$\bP^1$. From its definition, it follows that $\bD$ preserves any square 
root $\cK^{1/2}(\rho)$ of the canonical bundle of $\frM(x,\cB)$. Setting 
$h=-c, \mu=0$ gives a fixed-point for \eqref{dualine} and shows the 
Ansatz \eqref{factorline} to be correct. 
 
For more general admissible classes, $S$ and $\bD$ are only related
after splitting some filtrations. Denote by $\partial{E^*_x
V}/\partial H$ the sum of weight line bundles on $\frM(x,\cB)$ defined
by the virtual character $\partial\Tr_V/\partial H$ of $T$. We claim 
that
\begin{equation}\label{dualvect}
\begin{split}
\bD\,&\gr\left\{\cL(\mu+\rho)\otimes\exp[tE^*_\Sigma V]\right\} = \\
	&= \cL(S\mu +\rho)\otimes\gr\left
	\{\exp[-t\cdot \partial{E^*_x V}/\partial H)]\otimes
	\exp[tE^*_\Sigma V]\right\},
\end{split}
\end{equation}
for certain finite filtrations (term by term in $t$) on the two sides. 
To see this, let $\nu$ be the highest weight of $V$ and let $\cE'$ be 
the sheaf of sections of $E^*V$ whose $\lambda$-weight component vanishes 
at $x$ to order $\frac{1}{2}(\nu-\lambda)(H)$ or higher. This condition is 
stable under the $\alpha_0$-root $\mathfrak{sl}(2)$, so $\cE'$ descends 
to $\Sigma\times\frM'$. The quotient $\cQ=E^*V/\cE'$ is supported on 
$\{x\}\times\frM(x,\cB)$. It has a finite filtration whose associated 
graded sheaf is a sum of weight line bundles $\cO(\lambda)$ on $\frM(x,\cB)$, 
for the weights $\lambda$ of $V$ and various multiplicities. By construction,
\[
s_0(\gr\cQ) - \gr\cQ = \partial{E^*_x V}/\partial{H}.
\] 
Dualising $\gr\, E^*_\Sigma V$ along the fibres of $p$ then 
results in $\gr\, E^*_\Sigma V - \partial{E^*_x V}/\partial{H}$, proving 
\eqref{dualvect}.

The sign change of $\cI$ under the action \eqref{waffaction} of $S$ 
follows now by factoring the index map via $\frM'$, since splitting 
the filtrations does not change the index.
\end{proof}

\begin{example} When $\GC=\SL(2)$, $\frM(x,\cB)$ is the moduli stack of pairs 
$(\cE,L)$, where $\cE$ is a rank $2$ bundle with trivial determinant and $L$ 
a line in the fibre at $x$. $\frM'$ is naturally equivalent to
the stack of rank $2$ bundles with determinant identified with the
line bundle $\cO_\Sigma(-x)$. The morphism $p$ takes $(\cE,L)$ to its
sub-sheaf $\cE'$ of sections whose value at $x$ lies in $L$. The lines
$L$ assemble to the weight line bundle $\cO(-1)$ over $\frM(x,B)$. The
vector bundles associated to irreducible representations of $\SL(2)$ are 
the symmetric powers of $\cE$, and the maximal sub-sheaves in the proof 
of \eqref{waff} are the symmetric powers of $\cE'$. The quotient 
$S^n\cE/S^n\cE'$ is supported at $x$, its associated graded sheaf over 
$\frM(x,B)$ is $\bigoplus_{0\le k \le n} \cO(n-2k)^{\oplus k}$ and the
anti-symmetrisation is $\bigoplus\cO(k)^ {\oplus k}$ ($k=n\mod{2}$, 
$|k|\le n$). 
\end{example}

\begin{corollary}\label{delta}
The series $\cI$ represents a Weyl anti-symmetric linear combination of 
$\delta$-functions on $T_k$. In particular, it is supported at regular 
points only.
\end{corollary}

\begin{proof}
Weyl anti-symmetry is clear. Assume first that $\GC$ is simply connected, 
so that $\Pi = \pi_1T$. At $t=0$, invariant functionals under the 
lattice $\Pi$ are spanned by $\delta$-functions supported on $F$. However, 
the $t$-deformed action is obtained from the one at $t=0$ by the change of 
coordinates \eqref{coord}; so the $\Pi$-invariant Fourier series are 
spanned by the $\delta$-functions at the regular $f_t$. 

In general, the $\waff$-symmetry of \S\ref{waffaction} can be enhanced by 
the action of the co-weights of the centre $Z(G)\subset G$: geometrically, 
 these central co-weights define elementary transformations on bundles 
which translate the components of $\frM$, and the multiplicative factor 
in the $\waff$-action corrects for the change in $\cL\otimes\exp
[tE^*_\Sigma V]$. The extended lattice is co-compact in $\frt_k$, so 
our index functional is a span of $\delta$-functions, as before.
\end{proof}

\subsection{Odd classes.} \label{oddserre}
The arguments of this section also apply to more general 
bundles $\bigotimes_k E_{C_k}^*V_k\otimes \cL\otimes\exp[tE^*_\Sigma V]$ 
which include odd factors $E^*_{C_k}V_k$ from (\ref{atbott}.ii): each 
$C_k$ can be moved to avoid the Hecke point $x$, and $E^*_{C_k}V_k$  
remains unchanged in the Serre duality step. Let us rephrase this 
observation to match our indirect definition~\ref{indexdef} of the index. 

The index is obtained by first pushing down $\boxtimes_k E^*V_k\otimes
\cL\otimes\exp[tE^*\Sigma V] \otimes E^*_xV_\mu$ to $\Sigma^n$, 
and then taking the index over $\times{C}_k$. Now, $\times{C}_k$ lies 
within $(\Sigma^\circ)^n$, with $\Sigma^\circ: = \Sigma\setminus\{x\}$, 
so we can restrict our bundle to $\frM(x,\cB) \times (\Sigma^\circ)^n$. 
In repeating the arguments, 
we note that each bundle $E^*V_k$ is in fact pulled back from $\frM'
\times (\Sigma^\circ)^n$: this is because the corresponding $\cE'$ 
used the proof of \eqref{waff} is part of an exhaustive filtration, 
once we allow arbitrary poles at $x$. Therefore, each factor survives 
Serre duality unchanged, leading to the same symmetry.  

\section{Abelianisation}
\label{5}
We now prove Theorems~\ref{general} and \ref{local}. For technical 
reasons, we must use the stack $\frM(x,\cB)$ of bundles decorated 
with a Borel structure at $x\in \Sigma$. Fix an admissible class $\cE$ 
and let 
\[
\cI_\cE:= \sum\nolimits_\mu \Ind\left(\frM(x,\cB);
	\cE(\mu+\rho)\right)\cdot e^{-\mu}.
\]
Taylor expansion in Corollary~\ref{delta} (and \S\ref{oddserre}, 
if $\cE$ contains odd classes) shows that $\cI_\cE$ is a distribution 
supported on a finite set of regular points in $T_k$. Let  $\widetilde
{\frM}_T$ denote the moduli of $T$-bundles trivialised at $x$ (a 
disjoint union of Jacobians), and $\Ind_T$ the $T$-equivariant index 
of a vector bundle lifted from ${\frM}_T$. Call $\nu_B$ the virtual 
normal bundle of the morphism $\frM_T\to \frM(x,\cB)$. 

\begin{proposition}\label{abel}
Over the regular part of $T_k$, 
$\cI_\cE = \Ind_T\left(\widetilde{\frM}_T;\cE(\rho)/\lambda_{-1}
	(\nu_B^\vee)\right)$ as distributions.
\end{proposition} 

\noindent For the proof, we need a preliminary calculation of the 
Euler complex of \S\ref{spectral}.

\subsection{The Euler complex.} \label{eulerinverse}
For any $\xi$ labelling a Shatz stratum, recall the complex
\[
\on{Eul}(\nu_\xi)^{-1}_+ = \Sym\, \left(R\pi_* E^*
	(\frp_\xi/\frg_\xi)[1]^\vee \oplus R\pi_* E^* (\frg/\frp_\xi)[1]
	\right) \otimes \det(R\pi_* E^* (\frg/\frp_\xi)[1]) [d_\xi]. 
\]
It splits by $\xi$-eigenvalue into bounded complexes with coherent 
cohomologies, and for index purposes we may perform $K$-theoretic 
cancellations. One such arises from Serre duality 
\begin{equation}\label{pgdual}
R\pi_* E^*\frp_\xi/\frg_\xi[1]^\vee =
	R\pi_* (E^*\frg/\frp_\xi \otimes{K}), 
\end{equation}
using the $\GC_\xi$-isomorphism $(\frg/\frp_\xi)^\vee = \frp_\xi/\frg_\xi$. 
Replacing the second complex, in $K$-theory, by $R\pi_* E^*\frg/
\frp_\xi \oplus(2g-2) E^*_x \frg/\frp_\xi$ simplifies the symmetric 
factor to 
\[
\left(\Sym\, E^*_x\frg/\frp_\xi\right)^{\otimes(2g-2)}.
\]
Similarly, the determinant can be rewritten as $\det^{g-1}
E^*_x\frg/\frp_\xi\otimes D_\Sigma(\frg/\frp_\xi)$. Using Serre 
duality as in \eqref{pgdual} again, we get $D_\Sigma(\frg/\frp_\xi)
\cong D_\Sigma(\frp_\xi/\frg_\xi)$. With $\Delta_\xi, \rho_\xi$ and 
$\cK_\xi$ denoting the $\GC_\xi$-counterparts of $\Delta,\rho,\cK$, 
and using the co-dimension formula $d_\xi = (g-1)\dim(\frg/\frp_\xi) + 
2(\rho-\rho_\xi)(\xi)$, we obtain the following $K$-theoretic 
replacement of the inverse Euler complex:
\begin{equation}\label{inverseuler}
\bE(\nu_\xi)_+^{-1}:= (-1)^{2(\rho-\rho_\xi)(\xi)} E^*_x 
	\left(\Delta_\xi/\Delta\right)_+^{2g-2}\otimes(\cK_\xi/\cK)^{1/2},
\end{equation}
remarkably, a line bundle twist of a geometric series of weight 
line bundles. The subscript denotes that $\xi$-negative modes 
are to be chosen for the Fourier expansion of the Weyl denominator. 

For simplicity, we have avoided parabolic structures; in the special 
case of Borel structure at a single point $x$, to be used in the 
proof below, \eqref{inverseuler} carries an additional factor of 
$e^{\rho_\xi -\rho} (\Delta_\xi/\Delta)_+$, from the flag varieties of 
$G$ and $G_\xi$. 

\begin{proof}[Proof of Proposition \ref{abel}.]
Stratify $\frM(x,\cB)$ using a generic polarisation (\S\ref
{shatzborel}) and express the index distribution $\cI_\cE$ as 
a sum of contributions $\cI_{\cE,\xi}$ from Shatz strata $\frM_\xi$ 
as in \eqref{notlocalisation}. Each $\cI_{\cE,\xi}$ is a formal 
Fourier series whose coefficients are the indices of 
$\cE_\xi \otimes \bE(\nu_\xi)_+^{-1}$ over the moduli stack 
$\frM_{G_\xi,\xi}^{\on{ss}}$ of semi-stable $G_\xi$-bundles 
of topological type $\xi$, with Borel reduction at $x$. 
We claim that
\begin{enumerate}\itemsep0ex
\item Each $\cI_{\cE,\xi}$ is a distribution on $T_k$.
\item Unless ${G_\xi}= T$, $\cI_{\cE,\xi}$ is supported
in the $\frg$-singular locus of $T_k$.
\item $\sum \cI_{\cE,\xi} = \cI_\cE$, convergent as a series 
of distributions.\footnote{Convergence as a formal Fourier 
series is clear, but we need distributional convergence for 
our argument.}
 \end{enumerate}
When $\frg_\xi=\frt$, $\bE(\nu_B)_+^{-1} = \lambda_{-1}
(\nu_B^\vee)^{-1}$ over the regular locus of $T_k$, proving our 
proposition subject to the three claims.

We now prove the claims. Since the polarisation is generic, 
the stack $\frM_{G_\xi,\xi}^{\on{ss}}$ is the quotient by $T$ of a 
smooth, quasi-projective variety on which the centre $\frz_\xi$ of 
$\frg_\xi$ acts trivially and $\frt/\frz_\xi$ acts freely.  If we 
ignore the automorphisms coming from the trivial action of $Z_\xi$, 
then $\frM_{G_\xi,\xi}^{\on{ss}}$ is a smooth, proper Deligne-Mumford 
stack. The Narasimhan-Mehta-Seshadri construction \cite{ms:pb} presents
the underlying orbifold as a (locally free) quotient of a compact 
manifold by a compact group: namely the quotient by $(T/Z_\xi)_k
$-conjugation of the manifold $M_\xi^*$ of flat unitary 
$\GC_\xi$-connections on $\Sigma\setminus\{x\}$ with a prescribed, 
regular value in $T_k$ of the monodromy at $x$.

The index of a vector bundle over $\frM_{G_\xi,\xi}^{\on{ss}}$ is 
the $Z_\xi$-invariant part of the index of its direct image to the 
orbifold. GAGA, applied to the coarse moduli space, allows us to 
use the holomorphic Euler characteristic instead. As in the manifold 
case, that can be identified with the index of a twisted Dolbeault 
operator (see e.g. Duistermaat \cite{du:he}). By Kawasaki \cite
[Example II]{kaw:ind}, the latter is the invariant part of the 
distributional index of a twisted Dolbeault operator on 
$M_{\xi}^*$, which is transversally elliptic for the 
$T_k$-action.\footnote{The reader may refer to \cite{ve:eq, pa:lo} 
for further discussion of these methods.} By \eqref{inverseuler}, 
expansion into Fourier modes equates the index series $\cI_{\cE,\xi}$ 
with the distributional Dolbeault index of $\cE_\xi\otimes 
(\cK_\xi/\cK)^{1/2}$ on $M_\xi^*$, multiplied by $(\Delta_\xi/
\Delta)_+^{2g-2}$. Since $Z_\xi$ acts trivially, the distributional 
index is in fact a Fourier polynomial along $Z_\xi$, and the 
$\xi$-negative choice of the Fourier expansion ensures convergence 
of the sum to a distribution. This proves claim (i). 

We derive claim (ii) from Atiyah's localisation theorem \cite[Thm.~4.6]
{a}, which asserts that the distributional index of a transversally 
elliptic operator is supported over the union of all stabiliser subgroups. 
Now, the freedom of $\pi_1G$ implies that all stabilisers of the $T_k
$-action on $M_{G_\xi,\xi}^*$ lie in the $\frg$-singular locus. 
Indeed, a result of Borel's ensures that the $\GR$-centraliser of any 
$\frg$-regular torus element is $T_k$ itself; but if all monodromies 
were in $T_k$, then the monodromy around $x$ would be trivial. 

Finally, for (iii), it suffices to fix $\xi$ and $G_\xi$ and show 
convergence of the sum of $\cI_{\cE,\xi+\gamma}$ over the co-weights 
$\gamma$ of $Z_\xi$. We'll also need to divide into co-sets of $W/W_\xi$, 
for the different expansions of the inverse Euler class. Compared 
with $\cI_{\cE,\xi}$, the Atiyah-Bott factors (i) and (ii) in  
$\cI_{\cE,\xi+\gamma}$ are unchanged, while each index bundle 
$E^*_\Sigma V$ acquires a summand $E_x^*(\partial\Tr_V/\partial
\gamma)$. This is a sum of weight spaces of $V$, with multiplicities 
linear in $\gamma$, and factors out of the index. Finally, a line 
bundle $\cL$ of $\cE$ gets shifted by the weight $\iota(\gamma)h$ of 
$T$ (see for instance \eqref{first} below), while $\bE(\nu_\xi)^{-1}_+$ 
acquires a factor of $e^{\iota(\gamma)c}$ from the canonical bundles. 
This gives a sum of the form
\[
\sum_\gamma \cI_{\cE,\xi+\gamma} = 
	\sum_{\gamma;j;\mu} p_{j,\mu}(\gamma)\cI'_{\cE,j,\mu}\cdot e^\mu 
	\cdot{e}^{\iota(\gamma)h'}\prod_{\alpha(\xi+\gamma) < 0} 
	(1-e^\alpha)^{2-2g},
\]
over finitely many values of $j,\mu$, with distributional Dolbeault 
indices $\cI'_{\cE,j,\mu}$ of vector bundles over $M^*_\xi$ 
and polynomials $p_{j,\mu}$ in $\gamma$. The distributions 
$\cI'_{\cE,j,\mu}$ are in fact Fourier polynomials along $Z_\xi$, so 
the negativity constraint on the roots $\alpha$ and negative-definiteness 
of $h'$ on the co-weight lattice assures distributional convergence  
after expanding out into a Fourier series.
\end{proof} 

\subsection{The Jacobian contributions.}
In preparation for the proof of Theorem \ref{general}, we now spell out 
the Riemann-Roch formula for the $T$-Jacobians.  The components of 
$\frM_T$ are labelled by the first Chern classes of $T$-bundles, 
valued in the co-weight lattice. Each component $\frM_T^{(\gamma)}$ 
factors as $J_\gamma\times BT$, so that the projection to $BT$ lifts 
the $T$-representation $\bC_\mu$ with weight $\mu$ to the line bundle 
$E^*_x\bC_\mu$, and each $J_\gamma$ is identified with the $T$-Jacobian 
$J:=J_0$ by an elementary transformation at $x$.  Call $\omega$ the 
positive integral generator of $H^2 (\Sigma)$ and $\Psi$ the duality 
tensor in $H_1(\Sigma)\otimes H^1(\Sigma)$. After the natural 
identifications 
\[
H^1(J)\cong
T^\vee J\cong H^1(\Sigma; \frt)^\vee\cong \frt^\vee\otimes H_1(\Sigma),
\] 
we have on $\frM_T^{(\gamma)}\times\Sigma$
\begin{equation}\label{chern}
c_1(E^*\bC_\mu) = \pi^*c_1(E^*_x\bC_\mu) + 
		\mu(\gamma)\cdot\omega +	\mi\mu\otimes\Psi.
\end{equation}
With the cup-product form $\eta\in\Lambda^2 H_1(\Sigma)$, we note the 
relation
\[
(\mu\otimes\Psi)^2 = -2\mu^{\otimes{2}}\otimes\eta\wedge\omega 
\in H^4(J\times\Sigma),
\]
where $\mu^{\otimes{2}}$ is the square in $\mathrm{Sym}^2\frt^\vee$.
We now use the equivariant Chern character to convert admissible 
$K$-classes over $\frM_T = J\times BT$ into cohomology classes on 
$J$ with coefficients in $R_T$. For instance, the Chern character 
$\on{Ch}(E^*_x\bC_\mu)$ becomes the group character $e^\mu$. 
Formula \eqref{chern} gives
\begin{align}\label{chernchar}
\on{Ch}(E^*\bC_\mu) &= e^\mu\left(1+\mu(\gamma)\cdot\omega\right)
				(1+\mi\mu\otimes\Psi + \mu^{\otimes{2}}\otimes
				\eta\wedge\omega), \nonumber\\
\on{Ch}\left(E^*_\Sigma\bC_\mu\right) &= e^\mu\left(\mu(\gamma) +
				\mu^{\otimes{2}}\otimes\eta\right)\nonumber\\	
\on{Ch}\left(E^*_C\bC_\mu\right)&= 
			e^\mu\cdot C\otimes \mi\mu\in H^1(J),\\	
\on{Ch}(D_\Sigma\bC_\mu) &= e^{-\mu\cdot\mu(\gamma)}
				\exp(-\mu^{\otimes{2}}\otimes\eta), \nonumber\\
\on{Ch}\left(\exp\left[tE^*_\Sigma\bC_\mu)\right]\right) &= \exp\left\{
		te^\mu\left[\mu(\gamma)+\mu^{\otimes{2}}\otimes\eta\right]\right\},
			\nonumber
\end{align}
whence we get on $BT\times J_\gamma$, for any $T$-representations $U,V$, 
the two formulae
\begin{eqnarray} \label{first}
\on{Ch}(D_\Sigma U) &=& 
	e^{\iota(\gamma)h}\cdot\exp(h\otimes\eta), \label{chdet}\\
\label{second} \on{Ch}\left(\exp[tE^*_\Sigma V]\right)(u) &=& 
	\exp\left\{t\left[\partial\Tr_V(u)/\partial\gamma
		+ \hess_V(u)\otimes\eta\right]\right\},
\end{eqnarray}	
with the metric $h=-\Tr_U$ on $\frt$ and the Hessian $2$-form $\hess_V(u)$ 
of $\Tr_V$ at $u\in T$.  

Finally, to find the Riemann-Roch expression for $\lambda_{-1}(\nu^\vee)$, 
we apply the argument of \S\ref{eulerinverse} over a $\gamma$-component 
of $\frM_T$, restricting to the regular points of $T_k$ (where the 
choice of expansion of the series is immaterial). We get from \eqref{first}
\begin{equation} \label{third}
\on{Ch}(\lambda_{-1}\nu^\vee)^{-1}=(-1)^{2\rho(\gamma)}\Delta^{2-2g}
	e^{c(\gamma)}\exp[c\otimes\eta].
\end{equation}

\begin{remark}\label{injtrans}
Note from \eqref{chdet} that $h$ can be recovered from $c_1
(D_\Sigma{U})$ when $G$ is a torus, and then for any $G$ by 
passing to the maximal torus.
\end{remark}

\begin{proof}[Proof of Theorem \ref{general}]
Summing over $\gamma$ the products of contributions in \eqref{chdet}, 
\eqref{second} and \eqref{third} gives the following answer on $T\times J$:
\begin{align}\label{sum}
\sum_{\gamma\in\pi_1T} &\on{Ch}\left(\cL \otimes
	\exp[tE^*_\Sigma V]\right)\wedge \on{Ch}(\lambda_{-1}\nu^\vee)^{-1} =
	\nonumber\\
	&= \sum_\gamma\frac{(-1)^{2\rho(\gamma)}}{\Delta(u)^{2g-2}}
	\left[u^{h'}\exp\left[td\,\Tr_V(u)\right]\right]
	^\gamma\cdot{\exp\left\{[h'+t\hess_V(u)]\otimes\eta\right\}}=
	\nonumber\\
	&= \delta_\rho\circ\chi'_t(u)\cdot\left.\exp\left\{[h'+t\hess_V(u)]
	\otimes\eta\right\}\right/\Delta(u)^{2g-2}.
\end{align}
Observe now that 
\begin{equation}\label{skewgauss}
\int_J \exp\left\{[h'+t\hess_V(u)]\otimes \eta\right\} 
	=|F|^g \det\nolimits^g\left[1+t\hess_V(u)^\dagger\right].
\end{equation}
At $t=0$, this follows because $|F|$ is the determinant of $h':\frt\to\frt^\vee$ 
(with volume form normalised by the respective lattices), and the polarisation 
$\eta$ on the $\GL(1)$ Jacobian is principal; while from $t=0$ the formula 
is clear in general. Theorem \eqref{main} now follows from \eqref{thetas} 
and \eqref{deformtheta}.

To prove Theorem \ref{general}, recall from \eqref{chernchar} the Chern 
characters $\mi\, d\,\Tr_V (u)\otimes [C] \in H^1(J)$ of odd classes 
$E^*_CV$. Including a monomial $\psi$ in these classes in the integral 
\eqref{skewgauss} multiplies it precisely by the $[\psi](u)$ defined by the 
contraction procedure in \S\ref{odd}.
\end{proof}

\begin{remark} \label{discon}
Summing over the relevant part of $\frM_T$ gives the correct answer
for each component of $\frM_G$ separately. Similarly, we can produce 
a formula for the index over the moduli of vector bundles with fixed 
but non-trivial determinant from the sum over appropriate Jacobians. 
However, torsion in $\pi_1$ brings in additional contributions from 
principal bundles under the normaliser of $T$ in $\GC$; see the closely 
related calculation in \cite{amw} for line bundles.
\end{remark}

\begin{proof}[Proof of Theorem \ref{local}.]
In $K^0(\frM_T)$, $\nu_B= \nu + E^*_x(\frg/\frb)$, cf.~\S\ref{abelsubsect}, 
so $\lambda_{-1}(\nu_B) = \lambda_{-1}(\nu) \otimes \lambda_{-1}
(E^*_x(\frg/\frb))$. Anti-symmetry \eqref{waff} allows us to sign-average 
over $W$: Weyl's character formula converts $e^\mu/\lambda_{-1}
(E^*_x(\frg/\frb)^\vee)$ into $E^*_xV_\mu/|W|$ while leaving the factor 
$\cE/\lambda_{-1}(\nu^\vee)$ unchanged. 
\end{proof}

\section{Witten's formulae from the large level limit}
\label{witten}
Assume now that the genus $g$ is $2$ or more. If $\frM$ were a compact
manifold of complex dimension $d=(g-1)\dim G$, Riemann-Roch would
enforce the behaviour
\begin{equation}\label{asymp}
\Ind\left(\frM;\psi^n\cE\right) = n^d\int_{\frM} \on{Ch}(\cE) + O(n^{d-1})
\end{equation}
for any $K$-class $\cE$ and its $n$th Adams power $\psi^n \cE$. 
(Recall that $\psi^n L = L^n$ for a line bundle $L$, and $\psi^n$ 
extends to $K$-theory additively using the splitting principle.) 

In general, even the meaning of the integral on the right is unclear.
Suppose, however, that $\cE$ is a product of a polynomial in the Atiyah-Bott 
classes with a sufficiently large admissible line bundle. Then, for all 
$n$, $\Ind\left(\frM;\psi^n\cE\right)$ has vanishing contribution from 
the unstable strata (Lemma \ref{indexlemma} and \S\ref{vanish}), so 
the leading $n$ asymptotic term in the index comes from the semi-stable 
stratum in \eqref{notlocalisation}.  This contribution is slightly 
complicated by the singularities of the moduli space $M$. More precisely, 
the index of $\psi^n\cE$ over the semi-stable stratum is that of its 
direct image from $\frM^{\on{ss}}$ to $M$. Consider instead the direct image 
of the pull-back bundle $\tilde \cE$ to the orbifold desingularisation 
$\tilde M$ of $M$, obtained by Kirwan's method \cite{kir}. Because $M$ 
has rational singularities, the two indices agree, and when $\psi^n\cE$ 
descends to $\tilde M$, the leading term in the $M$-index is $n^d\int_
{\tilde M} \on{Ch}(\tilde\cE)$. Descent holds when all stabilisers 
on $\tilde{M}$ act trivially on the fibres, in particular $\psi^n\cE$ 
descends when all stabiliser orders in $\tilde{M}$ divide $n$.  

It is more convenient to find the leading term in the twisted 
limit $\cK^{1/2}\otimes \psi^n(\cK^{-1/2}\cE)$. Let $\cE = \cL\otimes 
\exp[tE^*_\Sigma V]$, specialising to even generators for simplicity.
Riemann-Roch implies
\[
\psi^n E^*_\Sigma V = \frac{1}{n} E^*_\Sigma(\psi^nV),
\]
and the properties of $\psi$ give
\begin{equation} \label{psin}
\psi^n\exp\left[tE^*_\Sigma V\right] = \exp\left[tE^*_\Sigma(\psi^nV)/n\right].
\end{equation}
Since $h'$ scales by $n$ and $d\,\Tr_{\psi^nV}(u) = n\cdot d\,\Tr_V(u^n)$, 
the transformation \eqref{coord} is unchanged, and the effect of the 
twisted $\psi^n$ operation is to pre-compose the map $\chi_t': T\to 
T^\vee$ in \eqref{deformtheta} with the $n$th power map on $T$. The key 
observation now is that the $n^d$ contribution to the sum in Theorem~\ref
{main}, as $n\to\infty$, come from those points $f_t$ located near the 
centre of $\GC$. Now, the descent condition on $\cE$ requires the centre 
of $G$ to act trivially on the fibres, with the result that the contributions 
near the various central elements agree, and summation over the centre 
can be concealed in the answer. Rescaling the $\log f_t$'s in the Lie algebra 
by $n$ recovers Witten's sum over integral weights in \cite[\S5]{wit}, with 
potential $Q = h'(\phi,\phi) + t\cdot\Tr_V(e^\phi)$ $(\phi\in \frg)$. For 
example, the rescaled Weyl denominator in the $\theta_t$ converges to 
the dimension formula for the representations. We only spell out the 
complete details for $\GC=\SL(2)$, but the method works in general 
(Remark \ref{genestimates}). 

Let $\cE$ be as above, with $c_1(\cL)=l\in H^2(M,\bZ)$ and $V$ of even 
spin $2j$. In the notation of \S\ref{su2}, a solution $\zeta_t$ of
\[
\zeta_t^{(2l+4)n}\cdot\exp\left[t\dot\varphi(\zeta_t^n)\right] = 1
\]
can be written
\begin{equation}\label{zeta}
\zeta_t = \exp\frac{\pi\mi k_t}{(l+2)n},
\end{equation}
where for each $k\in\bZ^+$, $k_t = k + k_1t+ k_2t^2 + \cdots$ formally 
solves the equation
\begin{equation}\label{kequ}
k_t + t\dot\varphi\left(\exp\frac{\pi\mi k_t}{l+2}\right) = k. 
\end{equation} 

\begin{proposition}\label{asymchern}
With $\tilde\cE$ as above, we have
\[
\int_{\tilde{M}} \on{Ch}(\cK^{-1/2}\otimes\tilde\cE) 
	= 2(l+2)^d\cdot\sum_{k=1}^\infty \left[1 + \frac{t\,\ddot
		\varphi(\exp\frac{\pi\mi k_t}{l+2})}{2l+4}
		\right]^{g-1}\cdot (\sqrt{2}\pi k_t)^{2-2g} \pmod{t^{(l+2)/j}}.
\]
\end{proposition}

\begin{remark}
\begin{enumerate}\itemsep0ex
\item Note that $l+2= c_1(\cK^{-1/2}\cE)$. 
\item To finite order in $t$, our formula involves integrals of 
polynomials in $\exp(c_1)$ and Atiyah-Bott cohomology classes; so our 
ingredients are equivalent to Witten's, the exponential term $\Tr(e^\phi)$ 
notwithstanding. But our truncation is needed precisely because the presence
of exponentials; for no $l$ does the formula hold to all orders in $t$.

\item Our answer seems at first to differ from \cite{wit}: the dimensions 
$k$ of the irreducible representations of $\SU(2)$ have been deformed to 
$k_t$. To reconcile the formulae, note that our first factor in the sum 
\eqref{asymchern} is the Jacobian determinant of the map $\xi\mapsto \xi 
+ t\nabla \Tr_V(e^\xi)$ on $\frt$, whereas its counterpart in \cite{wit} 
is the corresponding Jacobian on $\frg$. The ratio of the two is the 
volume ratio $k_t^2/k^2$ of the two co-adjoint orbits. 
\end{enumerate}
\end{remark}

\begin{proof}
We have $t$-truncated the formula to the place where unstable strata begin 
contributing to the index. We must then only check that \eqref{asymchern} 
gives the limiting $n^d$-coefficient in the index over $\frM$. To do so, 
we subdivide the summation range $1\le k < n(l+2)$ into an interior region 
and two ends. We then check that the interior sum is bounded by $o(n^d)$, 
while the ends gives the wanted $n^d$-contribution. 

The periodicity $k_t\mapsto k_t+2l+4$ shows that \eqref{kequ} involves 
only finitely many equations. All have analytic solutions 
for small $t$, so we can find a bound, independent of $n$ and $k$, for 
the variation of $k_t$ with $t$. This will allow us to replace $k_t$ 
by $k$ in some estimates. 

We now cut off at $k_- = \sqrt{n}$ and $ k_+ = (l+2)n - \sqrt{n}$.  
In-between, $|\zeta-\zeta^{-1}|> \pi/\sqrt{n}(l+2)$, so $|\theta^{1-g}| 
= O(n^{2g-2})$ and the sum over $k$ is bounded by $O(n^{2g-1})$, less 
than $o(n^{3g-3})$ when $g> 2.$
  
On the other hand, for $k< k_-$, Taylor expansion of the Weyl denominator 
gives
\[
\zeta_t - \zeta_t^{-1} = \frac{2\pi\mi k_t}{n(l+2)}\left(1
	+ O(n^{-1})\right),
\]
with $k$-independent error bound, so the $k$th term in the index sum is
\[
(l+2)^{2(g-1)}n^{3(g-1)}\cdot
	\left[\frac{2l+4 + t\,\ddot
		\varphi(\exp\frac{\pi\mi k_t}{l+2})}
		{(2\pi k_t)^2}\right]^{g-1}\left(1 + O(n^{-1})\right)
\]
and convergence of the series allows us to ignore the error. This and 
\eqref{asymp} give half of \eqref{asymchern}, the other half coming from 
the neighbourhood of $\zeta= -1$, by the central symmetry. 
\end{proof}

\begin{remark} 
The central symmetry relies on our choice of $V$ with even spin. 
Its absence for odd spin reflects the fact that the central automorphism 
of $\SL(2)$-bundles obstructs the descent of $E^*_\Sigma{V}$ to $\tilde{M}$. 
Similarly, there is an integration formula for the moduli of bundles with 
fixed determinant of degree $1$, which introduces a sign $(-1)^k$ in the 
sum  (cf.\ Remark~\ref{discon}). The level $l$ must now be even; 
else, the contributions near $\zeta=1$ and $\zeta=-1$ cancel instead 
of agreeing, even for line bundles. This reflects the fact that odd-level 
line bundles do not descend to the moduli space (again, the central 
automorphism acts by a sign).
\end{remark}

\begin{remark}\label{genestimates}
This argument works for any simple $G$. Subdivide the simplex $T_k/W$ 
of conjugacy classes into thickenings of width $1/\sqrt{n}$ of the faces. 
(First thicken the vertices, then the remainder of the edges, etc.) 
For each face $\Phi$, the factors in the Weyl denominator are bounded as 
above: $|\sin\alpha|> 1/n$ if the root $\alpha$ vanishes on $\Phi$, 
otherwise $|\sin\alpha|>1/\sqrt{n}$. With $Z_\Phi$ denoting the centraliser 
of $\Phi$, there are $\frac{1}{2}(\dim{Z}_\Phi -\dim{T})$ of the former 
and $\frac{1}{2}(\dim{G} - \dim{Z}_\Phi)$ of the latter. As $|F| =
O(n^{\dim{T}})$, the contribution of each point near $\Phi$ to the
index formula can be overestimated by $O(n^{(g-1)p})$, with
\[
p=\dim T + (\dim{Z}_\Phi -\dim T) + \frac{\dim{G}- \dim{Z}_\Phi}{2}
= \frac{\dim{Z}_\Phi + \dim{G}}{2}. 
\]
Even after adding $\dim \Phi$ to account for the number of terms, 
$(g-1)p$ is less than the dimension $(g-1)\dim{G}$ of $M$, unless 
$\Phi$ is a central vertex of $T_k/W$, so that $Z_\Phi = G$. Thus, 
only the $f_t$ near the centre contribute. The error estimate in 
their contribution proceeds as before. 
\end{remark}

\section{K\"ahler differentials}  
\label{kaehler}
In this section, we include the K\"ahler differentials $\Omega^\bullet$ 
over $\frM$ in our index. Recall that $\Omega^p = \lambda^p\left(R\pi_*
(E^*\frg \otimes K)\right)$, where $R\pi_*(\dots)$ is the (perfect) 
cotangent complex of $\frM$. Thus, $\Omega^p$ does not quite land  
in the admissible $K$-theory ring, but rather in its enlargement by 
the $\lambda$-operations. While the Abelian reduction formula \eqref
{local} and its proof carry over to this more general setting, our 
explicit index formula \eqref{general} does not immediately provide 
an answer. In this section, we show how to extend the formula  
to these more general $K$-classes.

As we will transfer the result to the moduli space $M$, we note 
the following improvement of Lemma~\ref{indexlemma}: the indexes of 
$\Omega^\bullet\otimes\cL\otimes \cE$ over the stack $\frM$ and over 
its semi-stable part $\frM^{\on{ss}}$ agree for large enough $\cL$, 
depending on the Atiyah-Bott monomial $\cE$ but \textit{not on the degree 
of the differentials}. The proof requires the finer calculation in \cite[\S7]
{tel2}. Note also that, for semi-simple $G$, the differentials on the 
stack of \textit{stable} bundles are the orbifold differentials over the 
moduli space of the same; in the reductive case, infinitesimal automorphisms 
cause a discrepancy which we leave in the care of the reader. 

Recall the notations of \S\ref{3}, in particular fix a representation 
$V$ of $G$. As $1+te^\alpha$ is a function on $T$, $(1+te^\alpha)^ 
\alpha$ is a $T^\vee$-valued map. Set
\begin{equation}\label{chi}
\chi'_{s,t}= \chi'\cdot e^{s\cdot d\,\Tr_V(.)}\cdot 
	\prod_{\alpha>0}\left[\frac{1+te^\alpha}
	{1+te^{-\alpha}}\right]^\alpha: T\to T^\vee. 
\end{equation} 
Denote by $F_{s,t}$ the set of solutions of the equation
\begin{equation}\label{points}
\chi'_{s,t}(f) = (-1)^{2\rho} \in T^\vee
\end{equation}
and by $F^{\mathrm{reg}}_{s,t}$ the subset of those which are 
\textit{regular} as $\GC$-conjugacy classes at $s=t=0$.  Call $H(f)$ 
the differential of $\chi'_{s,t}$ at $f\in T$; the notation $H$ stems 
from its agreement with the Hessian of the function on $\frt$
\[
\xi \mapsto \frac{h+c}{2}(\xi,\xi) + s\Tr_V(e^\xi) - 
\Tr_\frg\left(\mathrm{Li}_2(te^\xi)\right), 
\]
with Euler's dilogarithm $\mathrm{Li}_2$. Using the metric $(h+c)$, we 
convert $H$ to an endomorphism $H^\dagger$ of $\frt$ and define
\begin{equation}\label{theta}
\theta_{s,t}(f)^{-1} = |F|\cdot\prod_\alpha 
		\frac{1+te^\alpha}{1-e^\alpha}\cdot\det H^\dagger(f), 
\end{equation}
the product ranging over all roots. Note that $\det H^\dagger =1$ at 
$s=t=0$.

\begin{theorem}\label{index} With $\Omega_t:=\bigoplus_p t^p\cdot\Omega^p$, 
we have the index formula
\[
\Ind\left(\frM; \Omega_t\otimes\cL\otimes
		\exp[s E^*_\Sigma V]\otimes E^*_x U\right) = 
		(1+t)^{(g-1)\ell}\sum\nolimits_f \theta_{s,t}(f)^{1-g}\cdot \Tr_U(f),
\]
with $f\in F^{\mathrm{reg}}_{s,t}$ ranging over a complete set of Weyl 
orbit representatives. 
\end{theorem}

\begin{proof}
In topological $K$-theory, $R\pi_*(E^*\frg \otimes K) = E^*_\Sigma\frg
\oplus E^*_x \frg$, so $\Omega_t = \lambda_t\left(E^*_\Sigma\frg\right)
\otimes \lambda_t\left(E^*_x \frg\right)^{\otimes(g-1)}$. 
In terms of the Adams operations~$\psi^p$,  
\[
\lambda_t = \textstyle\sum_p t^p\lambda^p = 
\exp\left[-\textstyle\sum_{p>0} (-t)^p\psi^p/p\right].
\]
Using \eqref{psin}, we see that Theorem \ref{index} refers to the index 
over $\frM$ of
\[
\cL\otimes \exp\left[sE^*_\Sigma V - 
	\sum\nolimits_{p>0} \frac{(-t)^p}{p^2} 
		E^*_\Sigma\,\psi^p(\frg)\right]
		\otimes E^*_x\lambda_t(\frg)^{\otimes(g-1)}
		\otimes E^*_xU,
\]
which now has the form covered in Theorem \ref{main}. The associated 
equation
\[
\exp\left[ (h+c) + s\cdot d\Tr_V - \sum\nolimits_{p>0;\alpha} 
	\frac{(-t)^p}{p} e^{p\alpha}\cdot\alpha\right] = (-1)^{2\rho}
\]
is precisely \eqref{points}. To reduce formula \eqref{index} to 
Theorem~\ref{main}, observe that the pre-factor $(1+t)^\ell$ 
and the factors $1+te^\alpha$ in \eqref{theta} come from the character 
of $\lambda_t(\frg)$, which factors as $(1+t)^\ell\cdot\prod_\alpha 
(1+t e^\alpha)$. 
\end{proof}

\begin{remark} \label{skew}
Odd generators are included as in Theorem~\ref{general}, using the 
contraction procedure with the inverse of the bilinear forms $H(f)$.
\end{remark}

\subsection{Full-flag parabolic structures.}\label{fullflag}
A formula for the stack $\frM(x,\cB)$ follows by considering the 
projection $\frM(x,\cB)\to\frM$, with fibres $\GC/B$: 
replace $\Tr_U$ in \eqref{index} by 
\[
\Tr_U\cdot\prod_{\alpha>0}\frac{(1+te^\alpha)}{(1-e^\alpha)},
\]
and sum over all points of $F^{\mathrm{reg}}$ instead of Weyl orbits. 
The numerator accounts for the differentials on the fibres $\GC/B$, 
while the denominator and summation over $W$ together constitute the 
Weyl character formula.

\section{The Newstead-Ramanan conjecture}
\label{8}
In important special cases, all semi-stable bundles over $\Sigma$ are 
stable and then $M$ is a compact orbifold. This happens when $\GC=\GL(n)$, 
for the components of degree prime to $n$, or else if we enrich 
the bundle with a sufficiently generic parabolic structure.\footnote
{Note that the moduli of stable vector bundles of degree $d$ is also 
that of stable vector bundles of degree $0$ but with parabolic 
structure defined by the vertex $\mathrm{diag}[2\pi\mi d/n]$ 
of the Weyl alcove of $\mathfrak{gl}(n)$; cf.~ Example~\ref{fibres}.}  
Henceforth, we place ourselves in one of these favourable situations.  
Let $\ell^{\on{ss}}$ and $\ell^{\on{c}}$ be the semi-simple and
central ranks of $\GC$.  The following result generalises an old
conjecture of Newstead and Ramanan \cite{new, ram}.

\begin{theorem}
The top $(g-1)\ell^{\on{ss}} + g\ell^{\on{c}} $ rational Chern classes of $M$
vanish.
\end{theorem}

\noindent 
For rational cohomology, we can pass to finite covers with impunity 
\cite[\S7]{ab} and split $\GC$ as a product of a torus and simple groups; 
so the only content of the theorem concerns $\ell^{\on{ss}}$.  We will 
prove an equivalent result in topological $K$-theory. Let $\GC$ be 
semi-simple of rank $\ell$.

\begin{theorem} \label{Groth}
The top $(g-1)\ell$ rational Grothendieck $\gamma$-classes of $M$ vanish. 
\end{theorem}
  
\noindent The $\gamma$-classes are recalled below, along with the 
equivalence of the two theorems above. In some cases, such as 
$\GC=\GL(n)$, $\SL(n)$ or $\Sp(n)$, $M$ is known to be free of homology 
torsion \cite{ab}, and we get an integral result. It seems to be unknown 
whether $K(M)$ is torsion-free for other (e.g.\ simply connected) groups.

To prove Theorem \ref{Groth}, we pair the total $\gamma$-class $\sum
t^k\gamma^k$ of $TM$ against any test class $\cE$ in $K^0(M)$ and
show that we obtain a polynomial in $t$ of degree no more than 
$\dim M-(g-1)\ell$. Since the index over the orbifold $M$ varies 
quasi-polynomially in the Chern classes of $\cE$, it suffices to check 
this behaviour when $\cE$ contains a large line bundle factor, which we 
will do using using the index formula \eqref{index}. 

This strategy is not new, cf.~Zagier \cite{zag} for $\SL(2)$, 
but the integration formulae over $M$ turned out to be unwieldy. Our 
index formula seems to be a better fit; the reason is the abelian 
localisation \eqref{local}. Indeed, over $\frM_T$, the tangent complex 
to $\frM$ has a trivial summand of rank predicted by the vanishing. 
The proof then consists in checking that nothing in the index formula 
spoils the vanishing that is already apparent. Still, the method has 
limits: thus, we were unable to decide whether the $\gamma$-classes 
vanish in \textit{algebraic} $K$-theory.

\subsection{The $\gamma$-classes.}

For a complex vector bundle $V$ of rank $r$ over a compact space $X$, 
define the classes $\gamma^p(V)\in K^0(X)$ as the coefficients of the 
following polynomial of degree $r$:
\[
\gamma_t(V) = \sum\nolimits_p t^p \gamma^p(V) := (1-t)^r \lambda_
{t/(1-t)}(V),
\]
with the total $\lambda$-class $\lambda_s(V) = \sum s^p\lambda^p(V)$, as 
before. Note that 
$$\gamma_t(V\oplus W) =\gamma_t(V)\cdot\gamma_t(W)$$ 
for vector bundles $V$ and $W$, while 
$$\gamma_t(L) = (1-t) + tL$$ 
for a line bundle $L$; these conditions determine $\gamma_t$ from 
the splitting principle.  Also, 
$$\gamma^1(L) = L -1, $$ 
the $K$-theory Euler class of the line bundle, and in this sense
$\gamma_t$ is the total $K$-theory Chern class. The next exercise is
included for the reader's convenience.

\begin{proposition}
The following assertions are equivalent.
\begin{trivlist}\itemsep0ex
\item (i) The top $d$ rational Chern classes of $V$ vanish.
\item (ii) The top $d$ rational $\gamma$-classes of $V$ vanish.
\item (iii) The polynomial $\lambda_t(V)\in K^0(X;\bQ)[t]$ is divisible 
by $(1+t)^d$.
\end{trivlist}
\end{proposition}

\noindent When $K^0(X;\bQ)$ satisfies Poincar\'e duality with respect to a 
map $\Ind: K^0(X) \to \bQ$, these conditions are equivalent to

\noindent (iv) For every $W\in K^0(X)$, $\Ind(\lambda_t(V)\cdot W)\in 
\bQ[t]$ vanishes to order $d$ or more at $t= -1$. 

\begin{proof}
Equivalence of (ii) and (iii) is clear from the inversion formula 
$\lambda_t = (1+t)^r\gamma_{t/(1+t)}$. Observe next that in the ring 
$R$ of symmetric power series in variables $x_1,\dots,x_r$ the 
ideal $(e_{r-d+1}, \dots, e_r)$ generated by the top $d$ elementary 
symmetric functions is the intersection of $R$ with the ideal 
$(x_{r-d+1},\dots,x_r)\in \bQ[[x_1,\dots,x_r]]$. The transformation 
$$x_i\mapsto y_i:= e^{x_i} -1$$ 
defines an automorphism of $\bQ[[x_1, \dots,x_r]]$ which preserves
$(x_{r-d+1},\dots,x_r)$. It follows that $(e_{r-d+1}, \dots, e_r)$
agrees with the ideal of the top $d$ elementary symmetric functions in
the $y_k$. Let $x_k$ be the Chern roots of $V$; then,
$$\on{Ch}\,\gamma_t(E) = \prod(1 + ty_i), $$ 
so the $\gamma$-classes are the elementary symmetric functions in 
$y$, and we conclude that $\text{(i)}\Leftrightarrow\text{(ii)}$.
\end{proof}

\subsection{Reduction to Borel structures.}  
\label{borelenough}
We now show that if Theorem \ref{Groth} holds for moduli
of bundles with Borel structures, then it holds for all parabolic
structures.  Let $\frM(x,\cP)$ denote the stack of bundles with a
$\cP$-parabolic structure at $x$ and call $\pi: \frM(x,\cB) \to 
\frM(x,\cP)$ the projection \eqref{fiber}. Over $\frM(x,\cB)$, we 
have a distinguished triangle of tangent complexes
\[
T_\pi\frM(x,\cB) \to T\frM(x,\cB) \to \pi^*T\frM(x,\cP) \to 
	T_\pi\frM(x,\cB)[1],
\]
leading to an equality in $K$-theory, 
\[
\lambda_t\left(T^\vee\frM(x,\cB)\right) = 
	\lambda_t\left(T^\vee_\pi\frM(x,\cB)\right)
	\otimes \pi^*\lambda_t(T^\vee\frM(x,\cP)).
\]
The fibres $\cP/\cB$ of $\pi$ are flag varieties; they are smooth and
proper, with cohomology of type $(p,p)$. Hodge decomposition gives 
\[
R\pi_*\left[\lambda_t\left(T^\vee\frM(x,\cB)\right)\otimes\pi^*\cE\right] = 
	\lambda_t(T^\vee\frM(x,\cP))\otimes\cE\cdot \sum (-t)^p b_{2p}(\cP/\cB)
\]
where the $b_{2p}$ are the Betti numbers. For $t=-1$, the last factor 
is positive and so it does not affect the vanishing order of the index. 

\subsection{Limit of the index as $t \to -1$.} 
In Theorem~\ref{index}, the desired factor $(1+t)^{(g-1)\ell}$ appears 
explicitly in the index formula, so to prove Theorem~\ref{Groth} we must 
check that no singularities in $\theta_{s,t}(f)^{1-g}$ or in $H(f_t)$  
(cf.\ Remark~\ref{skew}) reduce the order of vanishing at $t=-1$. To do so, 
we study the roots of \eqref{points}.  When $h>0$ and $t=s=0$, $\chi'$ is 
an isogeny and all roots are simple. The following lemma will ensure 
that they remain simple for all $t\in (-1,0]$ and small $s$.

\begin{lemma}\label{hessian}
If $h>c$, $s$ is small and $t\in (-1,0]$, the differential $H = 
d\chi'_{s,t}$ is non-degenerate on $T_k$.
\end{lemma}

\begin{proof}
With $\hess_V(f)$ denoting the Hessian of $\Tr_V$ at $f$, we have
\begin{equation}\label{log}
H = (h+c) + s \hess_V(f) + t \sum\nolimits_\alpha 
			\frac{e^\alpha}{1+te^\alpha}(f)\cdot\alpha^{\otimes 2}.
\end{equation}
Note that $\alpha^{\otimes 2}$ is negative semi-definite, $t\le 0$ and 
$\Re\frac{e^\alpha}{1+te^\alpha} \ge -1$ for $|e^\alpha| =1$. As 
$\sum_\alpha \alpha^{\otimes 2} = -2c$, $H$ is bounded below by 
$(h-c) + s\hess_V$.
\end{proof}

\noindent Skew-adjointness of $\chi'$ for $s=0$ then keeps the solutions 
in the \textit{compact} torus $T_k$ for small variations in the real 
time $t$, and thus for all times $t\in [-1,0]$. Non-degeneracy of $H$ also 
shows that the $s$-dependence in \eqref{points} can be solved order-by 
order for all $t\in (-1,0]$, and keeps $F_{s,t}$ (with $s$ as the formal 
variable) in a formal neighbourhood 
of $T_k$. We will now show
that $H$ remains regular at $t=-1$, so the solution can be perturbed
analytically in $s$ even there. As certain regular solutions \emph{do}
wander into the singular locus of $T$ as $t\to -1$, we need to control
this behaviour. Let $f_t= f_{0,t}$.

\begin{lemma}\label{converge}
Let $f_t\in F$ be regular at $t=0$ but singular at $t=-1$. For 
small $x = \sqrt{t+1}$, $f_t$ has a convergent expansion 
\[
f_t = f_{-1}\cdot\exp\left[\textstyle{\sum}_{k>0}\; x^k\xi_k\right].
\] 
Moreover, $\beta(\xi_1)\neq 0$ for any root $\beta$ such that 
$e^\beta(f_{-1}) = 1$.
\end{lemma}
\noindent Thus, the tangent line to $f_t$ at $f_{-1}$ is regular in 
the Lie algebra centraliser $\frz$ of $f_{-1}$. We obtain
\begin{equation}\label{limits}\begin{split}
\lim_{t\to -1}\frac{1+te^\alpha}{1-e^\alpha}(f_t) &= 
	1 \qquad\text{for all roots } \alpha \text{ of } \frg,\\
\lim_{t\to -1}\left(\frac{e^\beta}{1+te^\beta}(f_t) + 
	\frac{e^{-\beta}}{1+te^{-\beta}}(f_t)\right) &= 
	-1 - \frac{2}{\beta(\xi_1)^2}\quad\text{for roots $\beta$ of }\frz.
\end{split}\end{equation}
The limiting value $H(f_{-1})$ in \eqref{log} is then the positive 
definite form $h + s\hess_V + \sum_\beta \frac{\beta^{\otimes 2}}{\beta
(\xi_1)^2}$, summing over roots of $\frz$. This excludes unexpected  
singularities in the index formula.

\begin{proof}[Proof of (\ref{converge}).]
At $t=-1$, equation \ref{points} simplifies to 
\begin{equation}\label{at-1}
\exp[h + s\cdot d\Tr_V(.)]= 1;
\end{equation}
however, the cancellation involved conceals multiple solutions on the
singular locus in $T$. The latter partitions $T_k$ into alcoves that
are simply permuted by the Weyl group. We claim that each
\emph{singular} solution of \eqref{at-1} is a limit of at least one
solution behaving as in Lemma \ref{converge}.  If so, then by
Weyl symmetry there must be such a solution from each adjacent alcove.
Now, every \emph{regular} solution of \eqref{at-1} is also the limit
of a regular solution of \eqref{points}: this is because it is the
limit of \emph{some} solution, and Weyl symmetry plus Lemma~\ref{hessian} 
ensures that the points of $F$ that are singular at $t=0$ stay so until 
$t=-1$.  Finally, recall that the solutions of \eqref{at-1} in a \emph
{closed} Weyl alcove are in bijection with the regular solutions of 
\eqref{points} in that alcove. Our claim then accounts for the $t=-1$ 
limits of \emph{all} regular points of $F$ and proves Lemma~\ref
{converge}.

To prove the claim, it suffices to find a \emph{formal} solution $f_t$ as 
in the Lemma. \comment{why?} As $t$ converges faster than $f_t$ becomes
singular, the function $\chi'_{s,t}$ converges to \eqref{at-1}, so 
equation \eqref{points} is verified to zeroth order precisely 
when $f_{-1}$ solves \eqref{at-1}. To obtain the constraint on $\xi_1$, 
we differentiate in $x$:
\[
\frac{d}{dx}\chi'_{s,t}(f_t) = \iota(\xi')
	\left[(h+c) + s\hess_V(f_t)\right] +
		\sum_\alpha \frac{[2x + t\alpha(\xi')] e^\alpha}
			{1+te^\alpha}(f_t)\cdot\alpha,
\]
with $\xi = \sum_k \xi_k x^k$. The limit at $x=0$ is found from \eqref
{limits} and leads to 
\begin{equation}\label{linear}
	\iota(\xi_1)\left[h + s \hess_V(f_{-1})\right] =
	\sum\nolimits_\beta \frac{\beta}{\beta(\xi_1)},
\end{equation}
summed over the roots $\beta$ of $\frz$. Its solutions are 
the critical points of the function 
\[
\frt\ni \zeta \mapsto  \frac{1}{2}\left[h(\zeta,\zeta) + 
	s \hess_V(\zeta,\zeta)\right] -
		\textstyle\sum_\beta \log|\beta(\zeta)|.
\]  
This function is real-valued for $s=0$, blows up on each walls of the 
Weyl chamber of $\frz$ and is dominated by the quadratic term at large 
$\zeta$, so a minimum must exist \textit{inside the chamber}. Further, 
the Hessian 
\[
h + \sum\nolimits_\beta \frac{\beta^{\otimes 2}}{\beta(\zeta)^2}
\]
is positive-definite, so the minimum is non-degenerate and the $s$-perturbed
equation can also be solved for small $s$. Continuing  to higher order in 
$x$, we get a recursive family of equations for $k>1$
\begin{equation}\label{recursion}
\iota(\xi_k)\left(h + s\hess_V + \sum\nolimits_\beta \frac
	{\beta^{\otimes 2}}{\beta(\xi_1)^2}\right) = 
		(\text{expression in } \xi_j,\:\: j<k),
\end{equation}
solvable because of the same non-degeneracy. This proves 
our claim and thus the lemma.
\end{proof}

\appendix
\section{Background on $\frM$}
\label{stratif}

For the more analytically minded, the stack $\frM$ admits the 
\emph{Atiyah-Bott presentation} as a quotient of the space of 
$(0,1)$-connections by the group of complex gauge transformations; 
but its underlying \textit{algebraic} structure is essential for us. 
The algebraic geometry of the stack  was discussed in \cite{bl, ls} 
and further properties were developed in \cite{tel1, tel2}. In 
particular, $\frM$ is covered by quotients of smooth varieties by 
reductive groups; many general properties of sheaf cohomology follow, 
without the need of simplicial topos theory as in \cite{tel1}.  
In this appendix, we quickly review the variants of $\frM$ with 
parabolic structures and discuss the topological $K$-theory of 
$\frM$.

\subsection{Parabolic structures.} \label{parab}
Call $\cB$ the \textit{Iwahori subgroup} of the loop group $\GC\pz$, 
consisting of those formal Taylor loops whose value at $z=0$ lies in 
a fixed Borel subgroup $B$. For any subset $\Phi$ of simple affine 
roots, let $\cP_\Phi$ denote the \textit{standard parabolic subgroup} 
of $\GC\pz$ generated by $\cB$ and by the root $\SL_2$ subgroups from 
$\Phi$. If $\Psi \subset\Phi$, then the quotient $\cP_\Phi/\cP_\Psi$ 
is isomorphic (possibly non-canonically) to a homogeneous space for 
a subgroup of $\GC$. 

\begin{example}
\begin{enumerate}
\item $\cP_\emptyset = \cB$. More generally, if $\Phi$ consists of 
(linear) roots of $\frg$, then $\cP_{\Phi}$ is the subgroup of formal 
Taylor loops whose value at $z = 0$ lies in the parabolic subgroup 
$\cP_\Phi\cap \GC$ of $\GC$.  

\item If $\Phi = \{\alpha_0\}$ the non-linear simple root, $\cP_\Phi$ 
has Lie algebra $\mathrm{Lie}(\cB) \oplus z^{-1}\frg_{\hroot}$. We 
have $\cP_\Phi/\cB \cong\bP^1$. This parabolic subgroup appears in 
the proof of Proposition~\ref{waff}.
\end{enumerate}
\end{example}

For distinct $x_1,\ldots,x_n \in \Sigma$ and $\cP_1,\ldots, \cP_n$
standard parabolics, let $\frM({\bf x};{\bf \cP})$ denote the moduli
stack of $\GC$-bundles with quasi-parabolic structures at $x_1,\ldots,
x_n$. These are $G$-bundles over $\Sigma$ with singularities at the 
$x_i$, but with a reduction of the gauge group to $\cP_i$ near $x_i$. 
When $G$ is semi-simple, the uniformisation theorem \cite[Theorem 9.5]
{ls} shows that this is the quotient of the product of \textit{generalised 
flag varieties} $\GC\pz/\cP_i$ by the gauge group $G[\Sigma \setminus 
\{ x_1,\ldots, x_n \}]$ of the punctured curve (cf.\ also \cite[\S9]{tel2}). 
Let $\cP_1',\ldots,\cP_n'$ be standard parabolics contained in
$\cP_1,\ldots,\cP_n$ respectively.  The projections $\GC\pz/\cP_i \to
\GC\pz/\cP_i'$ induce a fibration
\begin{equation} \label{fiber}
 \frM({\bf x};{\bf \cP}') \to \frM({\bf x};{\bf \cP})
 \end{equation}
with fibres $\cP_1/\cP_1' \times \dots \times \cP_n/\cP_n'$.  

\begin{example} \label{fibres}
\begin{enumerate}
\item When each of the parabolics is defined by a subset of the 
linear roots of $\frg$, $\frM({\bf x}; \cP)$ is the stack of $\GC$ 
bundles on $\Sigma$ with reductions to the parabolic subgroups 
$P_1,\ldots,P_n$ over $x_1,\ldots,x_n$. In this case, the $G$-bundles 
have no singularities and $\frM({\bf x};\cP)$ admits a forgetful 
morphism to the moduli stack $\frM$ with fibre $\GC/P_1 \times 
\dots \times \GC/P_n$.

\item For $\GC = \on{GL}(n,\bC)$, every parabolic subgroup of 
$\GC\pz$ is conjugate to one defined by linear roots of $\frg$, 
so all parabolic bundles can be described as vector bundles with a 
choice of flags at the marked points.
\item If $\Phi= \{\alpha_0\}$, then $\frM(x,\cP_\Phi)$ fibres over 
$\frM(x,\cB)$ with fibre $\bP^1$.  
\end{enumerate}  
\end{example}

\subsection{Shatz stratification.} \label{shatzborel}
Each stack $\frM({\bf x};\cP)$ is equivalent to a stack $\frM_\Gamma$ 
of equivariant bundles on a suitable Galois cover $\tilde{\Sigma} \to 
\Sigma$ \cite[\S2.2]{tw}. The Shatz stratification of $\frM_\Gamma$ 
induces a stratification on $\frM({\bf x}; \cP)$. This depends on the 
choice of the cover, but the dependence can be reduced to a choice of of 
\emph{polarisation} on $\frM({\bf x};\cP)$. For a Borel structure at 
a single point $x$, this is equivalent to a choice of finite-order, 
regular conjugacy class $\GR$; to label the strata, we choose a lifting 
$u\in T_k$.  

The co-weights $\xi$ labelling the Shatz strata of $\frM$ in \S\ref
{shatz} have a geometric meaning: every stable bundle in 
$\frM_{\GC_\xi,\xi}^{\on{ss}}$ has a unique Hermitian connection with 
constant, $\GC_\xi$-central curvature $2\pi\mi\xi$. The construction 
above shows that every stable bundle in $\frM_{\GC_\xi,\xi}(x,\cB)^
{\on{ss}}$ has a Hermitian connection over $\Sigma\setminus\{x\}$ with 
constant central curvature $2\pi\mi\xi$ and holonomy $u$ at $x$. The 
central part of $u$ stems from the curvature, while the projection to 
$\mathrm{Ad} G_\xi$ comes from the global monodromy.  

The index of an admissible class $\cE$ over $\frM({\bf x},\cP)$ breaks 
up as before into a sum over strata. There is also an extra factor 
in the Euler complex, relating the flag varieties of $\GC$ and $\GC_\xi$. 
The key finiteness result \eqref{indexlemma} applies to this more 
general setting, but the vanishing of unstable local cohomologies 
requires the line bundle $\cL$ to match the choice of stratification; 
see \cite[\S9]{tel2}. 

\subsection{$K$-theory of $\frM$.} \label{KofM}
The homotopy type of the stack $\frM$ (which, by definition, is that 
of the geometric realisation of an underlying simplicial scheme) is 
that of the space of continuous maps from $\Sigma$ to $BG$. (This is GAGA 
plus the Atiyah-Bott construction of holomorphic bundles.) But it is 
more natural to assign to $\frM$ the equivariant homotopy type given 
by the conjugation action of $\GR\subset \GC$ on the space 
$C_*(\Sigma, BG)$ of continuous maps based at $x\in\Sigma$ to 
$BG$. This space is a principal fibration over a product of copies 
of $G$, with fibre the group $\Omega G$ of based loops in $G$ \cite{ab}. 
Then, $K^\bullet(\frM)$ is defined to be the $\GR$-equivariant 
$K$-theory of $C_*(\Sigma, BG)$. It is an inverse limit of finite 
modules over the representation ring $R_G$, taken over the finite parts 
of a $\GR$-cellular model of $C_*(\Sigma, BG)$. Similarly, $K^\bullet
\left(\frM(x,\cB)\right)= K^\bullet_{T_k}\left(C_*(\Sigma, BG)\right)$; 
it is a module over $K^\bullet(\frM)$ via the natural projection, and 
$K^\bullet(\frM)$ is a split summand.

Another description of $K^\bullet(\frM)$ arises by exhausting $\frM$ 
with open sub-stacks of finite type. Such sub-stacks are presentable 
as quotients of quasi-projective manifolds by linear algebraic groups, 
and their topological $K$-theory can be defined from continuous 
vector bundles that are equivariant under the maximal compact part 
of the acting group. (This can be shown to be independent of the 
quotient presentation.) If we use the finite, open unions $\frM_
{\leq \xi}$ of Shatz strata to exhaust $\frM$, the argument of Atiyah 
and Bott (see \cite{land} for the $K$-theory version) shows the 
surjectivity of the restriction maps between the $K^\bullet(\frM_
{\leq\xi})$ and leads to the description
\begin{equation}\label{invlim}
K^\bullet(\frM) = \lim\nolimits_\xi K^\bullet(\frM_{\leq \xi}), 
\qquad \gr\, K^\bullet(\frM) = \bigoplus\nolimits_\xi K^\bullet(\frM_\xi).
\end{equation}
The two constructions of $K(\frM)$ just described can be related by 
presenting $\frM$ as a quotient $\frM_*/G$ of the stack of $G$-bundles 
with a framing over $x$ modulo the action of $G$ on the fibre: $\frM_*$ 
can be presented as a quotient of a pro-variety with the homotopy type of 
$C_*(\Sigma, BG)$ by a pro-unipotent group. 

Comparison with the stack $\frM_T$ of $T$-bundles gives more information. 
Consider for simplicity $\frM(x,\cB)$. When $\pi_1G$ is free, the 
stabilisers of the $T_k$-action on the complement of $C_*(\Sigma,BT)$ 
in $C_*(\Sigma, BG)$ are contained in the singular locus. Consequently, 
after inverting the Weyl denominator $\Delta$ in the coefficients of 
$K$-theory, the restriction  $j^*: K\left(\frM(x,\cB)\right) \to 
K(\frM_T)$ becomes an isomorphism, compatible with the inverse limit 
\eqref{invlim}. Poincare duality on $\frM_T$ and our index formula show 
that 
\[
(j^*)^{-1} = (-1)^{2\rho}\cdot\cK^{1/2}\Delta^{2g-2}\cdot j_*,
\] 
with $j_*$ defined using the finite-dimensional stack structure. 
However, our index formula carries the additional information that
inverting $\Delta$ does not damage the index.
\begin{remark} 
Rationally, $C_*(\Sigma, BG)$ is a product $\Omega G\times G^{2g}$. 
The rational cohomology factors \cite{ab} as 
\begin{equation}\label{factor}
H^\bullet(\frM) \cong H^\bullet_{\GR}(\Omega G) \otimes_R 
	H^\bullet_{\GR}(G)^{\otimes_R 2g},
\end{equation}
with $R = H^\bullet(BG;\bQ)$. A similar factorisation follows for rational 
$K$-theory, with $R= \mathbb{Q}\otimes R_G$, by using Chern characters 
and fixed-point formulae. It is tempting to suggest that the analogue 
of \eqref{factor} holds for integral $K$-theory when $\pi_1G$ is free, 
but we only know how to prove this for the groups $\GL,\SL$ and $\Sp$. 
\end{remark}

\vskip.5cm
\small{
\noindent\textsc{C.\ Teleman}, \texttt{c.teleman@ed.ac.uk}\\
School of Mathematics, Mayfield Road, University of Edinburgh, 
Edinburgh EH9 3ED, UK. 
\\
\textsc{C.T.\ Woodward}, \texttt{ctw@math.rutgers.edu}\\
Dept.~of Mathematics, Rutgers University, 110 Frelinghuysen Road, 
Piscataway NJ 08854, USA. }
\end{document}